\documentclass[12pt, english]{smfart}
\usepackage{amsfonts, amssymb, amsmath, latexsym}
\usepackage{mathptmx}
\usepackage[english]{babel}
\usepackage{mysectionsen}
\usepackage[latin1]{inputenc}
\usepackage[T1]{fontenc}
\usepackage{textcomp}

\usepackage{graphicx}
\usepackage[dvips]{color}
\input xy
\xyoption{all}

\newcommand{\Int}{\mathrm{in}}
\newcommand{\codim}{\mathrm{codim}}

\theoremstyle{definition}
\newtheorem{thm}{Theorem}[section]

\newtheorem{rmk}[thm]{Remark}
\newtheorem{example}[thm]{Example}

\newtheorem{conve}[thm]{Convention}

\usepackage{verbatim}

\date{}

\DeclareMathOperator{\Gr}{Gr}

\DeclareMathOperator{\spa}{Span}
\DeclareMathOperator{\initial}{in}
\DeclareMathOperator{\rk}{Rank}
\DeclareMathOperator{\coker}{Coker}

\DeclareMathOperator{\id}{id}

\newcommand{\te}[1]{\text{#1}}

\begin{document}

\thispagestyle{empty}

\begin{center} {\Large \bfseries{\textsc{Images of Rational Maps of Projective Spaces.}}}

\end{center}

\vspace{0.3cm}

\begin{center}
\textsc{Binglin Li}
\end{center}

\vspace{0.3cm}

\begin{center} \today
\end{center}

\vspace{5cm}

\hrule

\vspace{0,5cm}

{\small
\noindent
{\bf Abstract:}
Consider a rational map from a projective space to a product of projective spaces, induced by a collection of linear projections. Motivated by the the theory of limit linear series and Abel-Jacobi maps, we study the basic properties of the closure of the image of the rational map using a combination of techniques of moduli functors and initial degenerations.  We first give a formula of multi-degree in terms of the dimensions of intersections of linear subspaces and then prove that it is Cohen-Macaulay. Finally, we compute its Hilbert polynomials.
\vspace{0,2cm}


\vspace{0,1cm}

}

\vspace{0,5cm}

\hrule

\newpage

\section{Introduction}
\subsection*{Main Result}
Let $V$ be a vector space over an algebraically closed field $k$, and consider the following rational map:
$$f:\ \  \mathbb P(V)\dashrightarrow \Pi_{i=1}^n \mathbb P(V/V_i),
$$
which is induced by the choice of linear subspaces $V_i\subset V$. Let $X(V_1,\cdots, V_n)$ be the closure of the image of $f$.  The purpose of this paper is  to study the basic properties of $X(V_1,\cdots, V_n)$, using techniques of moduli functors and initial degenerations.

Fix a choice of  $n$ linear subspaces $\{V_i\}_{i=1}^n$ with $\bigcap_{i=1}^n V_i=\{0\}$, and let $d_I$ be the dimension of $\cap_{i\in I}V_i$ and define $M(p)$ to be the set:
$$\{(m_1,\cdots, m_n)\in \mathbb Z^n_{\geq 0}: \ \ r+1-\sum_{i\in I}m_i> d_I \ \ \forall I\subset [n],\ \sum_{i=1}^n m_i=p\}$$

In this paper, we have:

\begin{Thm} Set $p=\max\{h; M(h)\neq\emptyset\}$. The dimension of  $X(V_1,\cdots, V_n)$ is $p$. Its multi-degree function takes value one at the integer vectors in $M(p)$ and zero otherwise. The Hilbert Polynomial of $X(V_1,\cdots, V_n)$ is
$$  \sum_{S\subset M(p)}(-1)^{|S|-1}\Pi_{i=1}^n \binom{u_i+\ell_{S,i}}{\ell_{S,i}}
$$
where $u_i$'s are the variables and $\ell_{S,i}$ is the smallest $i$-th component of all elements of $S$. Moreover, $X(V_1,\cdots, V_n)$ is Cohen-Macaulay.
\end{Thm}

\subsection*{Background and Motivation}

Questions concerning some basic geometric properties of $X(V_1,\cdots, V_n)$ naturally arise in many settings, some of which are detailed below.

First, in \cite{est-oss} by Esteves and Osserman,  the closures of rational maps of the following form
 $$ \mathbb P(V)\dashrightarrow \mathbb P(V/V_1)\times\mathbb P(V/V_2)
 $$
 arise as irreducible components of the closed subscheme of the fiber of the Abel map associated to a limit linear series. Moreover, the union of those irreducible components forms a flat degeneration of a projective space. In this work, they work with  singular curves of compact type with two components, so they only need to study the closure of the image in a product of two projective spaces, where $V_1\cap V_2=\{0\}$.

Secondly, in computer vision ~\cite{computer} by Aholt, Sturmfels and Thomas, an important theme is to study the geometry of the closure of the rational map:
$$  \mathbb P^3\dashrightarrow (\mathbb P^2)^n.
$$
In their case, all the $V_i$'s are one dimensional and have pairwise trivial intersection.

Moreover, in the theory of Mustafin varieties in \cite{mus} and \cite{non-arc-unif}, flat degenerations of $\mathbb P^r$ are carefully studied, and the spaces $X(V_1,\cdots, V_n)$ arise as irreducible components of the Mustafin degenerations.

Finally, $X(V_1,\cdots, V_n)$  manifest themselves in algebraic statistics in  \cite{sta} by Morton .

 Despite the ubiquity of the spaces $X(V_1,\cdots, V_n)$ , a systematic study of its the geometry (multi-degree, singularity type, initial degeneration, Hilbert Polynomial etc.) has not be conducted in full generality. Such generality is needed. The author's original motivation is to generalize the work of Esteves and Osserman to  arbitrary curves of compact type. In this case, one cannot just focus on the map to a product of two projective spaces, and also the vector subspaces $V_i$'s do not mutually just have trivial intersections anymore, but rather the only restriction we may impose is that $\bigcap _{i=1}^n V_i=\{0\}$.  We don't put restrictions on the dimension of $V$ and $V_i$'s.  Moreover, we hope that our generality can lead to better understanding of irreducible components of Mustafin Degenerations.

\vspace{.3cm}

\subsection*{Structure of the Paper}
In Section II, we will describe two closed subschemes $\mathcal Z_1$ and $\mathcal Z_2$ of $\Pi_{i=1}^n \mathbb P(V/V_i)$, and prove that they are equal and that they agree with $X(V_1,\cdots,V_n)$ set-theoretically.

In Section III we will compute the multi-degree of $\mathcal Z_1$. It seems to be an elementary linear algebra argument, but the proof involves moduli-theoretic techniques to exhibit the existence of the system of linear subspaces satisfying the desired property and such technique is also used to compute the deformation.

In the Section IV we will compute the initial degeneration of $\mathcal Z_2$ under a prescribed term order. The slogan is that multi-degrees determine the initial degenerations of $X(V_1,\cdots, V_n)$.  We also prove that $\mathcal Z_1$, $\mathcal Z_2$ and $X(V_1,\cdots, V_n)$ are isomorphic as schemes, and finally conclude that $X(V_1,\cdots, V_n)$ is Cohen-Macaulay.

Section V will be devoted to computing the multi-variable Hilbert polynomials.  Our computation of the Hilbert polynomial is not by counting monomials using the information of the initial ideal obtained in Section III, but rather by directly computing the Hilbert polynomial of its initial degeneration, which is a union of product of projective spaces whose intersections are also products of projective spaces, which is easier and more intuitive to deal with.

\subsection*{Acknowledgements}  The author would like to thank his adviser Brian Osserman for his detailed and patient guidance.  The author also wants to thank Professor Bernd Sturmfels for helpful discussions relating this work with \cite{mus} and \cite{sta}, and also Professor Allen Knutson for introducing the author \cite{mulfree} which reduces problem of ``Cohen-Macaulay''ness to the proof of ``multiplicity free''.

Also special thanks to Naizhen Zhang and Christopher Westenberger for encouraging my studies in algebraic geometry and Federico Castillo for explaining to me some foundational material in combinatorial commutative algebra and useful discussions. The author also wants to thank Michael Gr\"{o}chenig for carefully proofreading this paper.

\section{descriptions of $X(V_1,\cdots, V_n)$}

\subsection{Set-theoretic Description of $X(V_1,\cdots, V_n)$}\label{initial model}

In this section, we first give a set theoretic description of $X(V_1,\cdots, V_n)$. A scheme theoretic description of $X(V_1,\cdots, V_n)$ will be given in Section 2.3.   Observing that $k$-valued points of  $\Pi_{i=1}^n\mathbb P(V/V_i)$ correspond to $n$- tuples  $(W_1,\cdots, W_n)$, where each $W_i$ is a linear subspace of $V$ which contains $V_i$ as codimension one linear subspace, we have the following theorem.
\begin{Prop}\label{closed points F_1}\label{set}The  closure of the image of the rational map $\mathbb P(V)\dashrightarrow \Pi_{i=1}^n\mathbb P(V/V_i)$ are set-theoretically in bijection with $n$-tuples $(W_1,\cdots,W_n)$ where $W_i$ contains $V_i$ as a codimension one subspace and for any $I\subset [n]$, we have $\bigcap_{I}W_i\supsetneqq \bigcap_{I}V_i$.
\end{Prop}

Before proving the statement, let's first start with an example:
\begin{example}
Let $V$ be a five-dimensional linear subspace over a field $k$ with basis $\{e_i\}_{i=1}^5$. Let $V_1=\text{span}\{e_1,\ e_2\},\ \ V_2=\text{span}\{e_1,\ e_3\}$, and $V_3=\text{span}\{e_5\}$. Let $W_1=V_1\oplus\text{span}\{e_3\},\ \ W_2=V_2\oplus\text{span}\{e_4\}$ and $W_3=V_3\oplus\text{span}\{e_1\}$.

 We have  $e_1\in  W_1\cap W_2\cap W_3$, but $e_1\in (V_1\cap V_2)\setminus V_3$, which means $(W_1,W_2,W_3)$ is not in the image of $f$. We also have $e_3\in W_1\cap W_2$, and $e_3\in V_2\setminus V_1$. Consider the following one-parameter family: $\mathcal W_t=(\spa(V_1,(e_3+te_4)), W_2, \spa(V_3,(e_1+te_3+t^2e_4))$. When $t=0$, it is the point $(W_1, W_2, W_3)$, when $t\neq 0$, $\mathcal W_t\in \text{Im}f$ except for possibly finite many values of  $t$. In this example, $(W_1, W_2,W_3)$ is not in the image of $f$, but is in the closure of the image of $f$ \\
\end{example}

\emph{\textbf{Proof of Proposition \ref{closed points F_1}:}}

Note that $k$-valued points in the image of the rational map $\mathbb P(V)\dashrightarrow \Pi_{i=1}^n\mathbb P(V/V_i)$ correspond to $(\spa (V_1, v),\cdots, \spa(V_n, v))$ where $v\in V\setminus (\cup_{i=1}^n V_i)$, thus one containment follows. For the other containment,  given $C_0=(W_1,\cdots, W_n)$ satisfying $\cap_{i\in I}W_i\supsetneqq \cap_{i\in I}V_i$, for all $I\subset[n]$, we will construct a one-parameter family $C_t$, such that when $t\neq 0$, we have $C_t=(\spa(V_1,v_t),\cdots, \spa(V_n,v_t))$ with $v_t\in V\setminus (\cup_{i=1}^nV_i)$ except for possibly finite many values of $t$, and $C_t=C_0$ when $t=0$.

Let $I_0=[n]=\{1,2\cdots,n\}$, then take a $ w_0\in \cap_{i\in [n]}W_i\setminus \cap_{i=1}^n V_i$. Let $I_1=\{i|\ w_0\in V_i\}$. Then take a $ w_1\in \bigcap_{I_1}W_i\setminus\bigcap_{I_1}V_i$, and let $I_2=\{i|\ w_1\in V_i, \  i\in I_1\}$, and take $w_2\in \bigcap_{i\in I_2}W_i\setminus(\bigcap_{i\in I_2}V_i)$. Repeating this process, i.e. choose a $w_j\in \bigcap_{I_j}W_i\setminus\bigcap_{I_j}V_i$, and let $I_{j+1}=\{i\in I_j|\ w_j\in V_i\}$ so that we have
 $$ I_1\supsetneqq I_2 \cdots\supsetneqq I_j\cdots
 $$
  Each $I_k$ is finite, thus for some $m$ we have $I_{m+1}=\emptyset$.  For $i\in I_{j}\setminus I_{j+1}$, we can express $W_i=\spa(V_i,w_j)$ (note that $i$ and $j$ are not necessarily the same).

Now consider the following set parameterized by $t\in k=\bar{k}$:

It suffices to consider the one parameter family component wise:  for component $i$, we can find $j$, such that $i\in I_j\setminus I_{j+1}$. Let $(W_i)_t= \spa(V_i,\sum_{j\leq \ell\leq m}t^{\ell-j}w_{\ell})$.

When $t=0$, we have  $C_0=(W_1,\cdots W_n)$  the point we start with.    When $t\neq 0$, we have $(W_i)_t=\spa(V_i,\sum_{j\leq \ell\leq m}t^{\ell-j}w_{\ell})=\spa(V_i,\sum_{0\leq \ell\leq m} t^{\ell}w_{\ell})$, thus $((W_1)_t,\cdots, (W_n)_t)$ is in the image of $f$ except for possibly finitely many values of $t$.
Then the proposition follows.
\hfill $\Box$

\subsection{Description II}\label{matrix}

Given coordinate systems for all $V$ and $V/V_i$'s, let $A_i$ be the matrix representation of the linear map $V\rightarrow V/V_i$, and $q_i$ be $(x_{i,1},\cdots, x_{i, r+1-d_i})^T$. Set $I=\{\delta_1, \cdots, \delta_{|I|}\}$ to be a subset of $[n]$,  and consider the following matrix associated to $I\subset [n]$,
$$B_I=\begin{bmatrix}
  A_{\delta_1}&q_{\delta_1}&0&\hdots&0 \\
  A_{\delta_2}&0&q_{\delta_2}&\vdots&0\\
  \vdots&\vdots& &&&\\
  A_{\delta_{|I|}}&0&\hdots&\hdots&q_{\delta_{|I|}}
\end{bmatrix}$$
Denote the ideal generated by all of the $r+1-d_I+|I|$-minors of $B_I$ as $\mathfrak I_{B_I}$, and denote the ideal generated by $\{\mathfrak I_{B_I}|I\subset [n]\}$ as $I_f$. Description II of $X(V_1,\cdots, V_n)$ is the closed subscheme of $\Pi_{i=1}^n \mathbb P(V/V_i)$ cut by $I_f$.

\begin{Def} \label{g_I|_k} Let $I$ be a finite index set and $\{V_i\}_{i\in I}$ be a system of linear subspaces of $V$, and $V\supset W_i\supsetneqq V_i$ with $\dim W_i-\dim V_i=1, \forall i\in I$. Then define $g_I$ to be the following map:
$$g_I:\ \ \ V\bigoplus (\oplus_{i\in I}W_i/V_i)\rightarrow \oplus_{i\in I}(V/V_i)
$$
 induced by the following data:
\begin{itemize}
  \item[(1)] $V\rightarrow V/V_i$ the natural quotient map.
  \item[(2)] $g_{i,i}:\ W_i\rightarrow V/V_i$ is the natural composed map $W_i\hookrightarrow V\rightarrow V/V_i$, and $g_{i,j}:\ W_i\rightarrow V/V_j$ is the zero map for $i\neq j$.
\end{itemize}
\end{Def}

Given any $W_i\subset V$ which contains $V_i$ as a codimension one subspace, $W_i/V_i$ corresponds to a rational point of $\mathbb P(V/V_i)$ through the natural map $W_i/V_i\rightarrow V/V_i$.  Conversely, given a rational point of $\mathbb P(V/V_i)$, which corresponds to a one-dimensional subspace $\ell$ of $V/V_i$, take the preimage of $\ell$ under the natural quotient map $V\rightarrow V/V_i$, then one get a linear subspace $W_i$ of $V$ containing $V_i$ as a codimension one subspace. Given a homogeneous coordinate of $\mathbb P(V/V_i)$, which is equivalent to giving a basis of $V/V_i$, the matrix representation of $s_i:\  W_i/V_i\rightarrow V/V_i$ is the same as the homogeneous coordinate of the point corresponding to $W_i$ up to a scalar. One can see that the matrix representation of $g_I$ is the same as $B_I$.

One can immediately get that Description  II agrees with $X(V_1,\cdots, V_2)$  at least set-theoretically:
\begin{Prop}\label{k-value} Let $V$ be a vector space over a field $K$ with dimension $(r+1)$. Let $I$ be a finite index set and $\{V_i\}_{i\in I}$ be a system of linear subspaces of $V$, and $V\supset W_i\supsetneqq V_i$ with $\dim W_i-\dim V_i=1, \forall i\in I$. Then $\cap_{i\in I}W_i\supsetneqq \cap_{i\in I}V_i$ if and only if the following naturally induced map:
$$g_I:\ \ \ V\bigoplus (\oplus_{i\in I}W_i/V_i)\rightarrow \oplus_{i\in I}(V/V_i)
$$
has  $\bigwedge^{(r+1-\dim(\cap_{i\in I}V_i)+|I|)}g_I=0$, where $g_I$ is defined in Definition \ref{g_I|_k}.

\end{Prop}
\emph{\textbf{Proof}}

The condition $\bigwedge^{(r+1-\dim(\cap_{i\in I}V_i)+|I|)}g_I=0$ is equivalent to the condition that $g_I$ has nontrivial kernel. The kernel of $g_I$ is $(\cap_{i\in I}W_i)/(\cap_{i\in I}V_i)$.  The kernel is nontrivial if and only if $(\cap_{i\in I}W_i)/(\cap_{i\in I}V_i)\neq {0}$, which is equivalent to the condition that $\cap_{i\in I}V_i\supsetneqq \cap_{i\in I}W_i$.

\hfill $\Box$

\vspace{.3cm}

\subsection{Equivalence of Two Descriptions }\label{comparison}

In this section I will prove that not just set-level, but also as schemes, the two descriptions in subsection \ref{initial model}
 and \ref{matrix} are the same. What we have done essentially in Proposition \ref{k-value} is compare $k$-valued points of the two descriptions. In this section, in order to compare the scheme structures of the two descriptions, we are going to compare their functor of points.

 Before giving the definition of $\mathfrak F_1$ and $\mathfrak F_2$, we first define a moduli functor $\mathfrak F$, which contains $\mathfrak F_1$ and $\mathfrak F_2$ as closed sub-functors.
 \begin{Def} Define $\mathfrak F$ to be a functor :
 $$\mathbf{(Sch/k)}^{\mathrm{op}}\rightarrow \mathbf{Sets}
 $$
 For any scheme $S/k$, we associate the set  of tuples sub-vector bundles $(\mathcal W_1,\dots,\mathcal W_n)$ of $\mathcal V_S:=V\otimes \mathcal O_S$ with $\rk (\mathcal W_i)=\dim V_i+1$, for any $i\in [n]$. And for each $i\in [n]$, the sub-bundle $\mathcal W_i$ contains $V_i\otimes\mathcal O_S$ as a sub-bundle.
 \end{Def}
 Note that the moduli functor $\mathfrak F$ represents nothing else than $\Pi_{i=1}^n \mathbb P(V/V_i)$. The following two functors $\mathfrak F_1$ and $\mathfrak F_2$ represent closed subschemes of $\Pi_{i=1}^n \mathbb P(V/V_i)$.

 \begin{Def} Define $\mathfrak{F_1}$ to be a functor :
 $$\mathbf{(Sch/k)}^{\mathrm{op}}\rightarrow \mathbf{Sets}.
 $$
 Given a scheme $S/k$, we associate the set of tuples sub-vector bundles $(\mathcal W_1,\dots,\mathcal W_n)$ of $\mathcal V_S:=V\otimes \mathcal O_S$ with $\rk (\mathcal W_i)=\dim V_i+1$ for any $i\in [n]$ such that:
 \begin{itemize}
 \item[(i)] for each $i\in [n]$, the sub-bundle $\mathcal W_i$ also contains $V_i\otimes\mathcal O_S$ as a sub-bundle
 \item[(ii)] for any $I\subset [n]$ and the following bundle map:
    $$f_I:\ \ \ V/(\cap_{i\in I}V_i)\otimes \mathcal O_S\rightarrow \oplus_{i\in I}\mathcal V/\mathcal W_i,
    $$
    we have
    $$\bigwedge^{(r+1-d_I)}f_I=0$$
 \end{itemize}

\end{Def}

\vspace{.3cm}

In order the define the second functor, we need the following lemma:

\begin{Lemma}\label{functor two definition} Let $V$ be a vector space over $k$ and $V'\subset V$ is a sub vector space of $V$. For a scheme $S/k$ , consider the following exact sequence:
$$\xymatrix{0\rightarrow V'\otimes\mathcal O_S\ar[r]^-{i_1}&\mathcal W\ar[r]^-{i_2}&V\otimes\mathcal O_S},
$$
where $\mathcal W$ is a sub-bundle of $V\otimes \mathcal O_S$ with $\rk(\mathcal W)=\dim V'+1$, and all of $i_1,\ i_2$ and $i_2\circ i_1$ have locally free cokernels. Denote  $\mathcal L$ to be the cokernel of $i_1$, which is a line bundle. Then we have
 \begin{itemize}
 \item[(i)] the induced morphism $0\rightarrow \mathcal W\rightarrow(V/V')\otimes \mathcal O_S$ factors through $\mathcal L\rightarrow (V/V')\otimes\mathcal O_S$
 \item[(ii)] this morphism $\mathcal L\rightarrow (V/V')\otimes\mathcal O_S$ is a sub-bundle morphism, i.e., has locally free cokernel
 \end{itemize}
\end{Lemma}

\noindent \emph{\textbf{Proof:}}

For (i), consider the sequence $0\rightarrow V'\otimes\mathcal O_S \rightarrow \mathcal W\rightarrow V\otimes\mathcal O_S\rightarrow (V/V')\otimes \mathcal O_S$, then we get $0\rightarrow V'\otimes\mathcal O_S \rightarrow \mathcal W \rightarrow (V/V')\otimes \mathcal O_S$ in which the morphism from $V\otimes \mathcal O_S$ to $(V/V')\otimes\mathcal O_S$ is zero.  Therefore by the universal property of cokernel, we naturally get the following commutative diagram:
$$\xymatrix{ \mathcal W\ar[rr]^-{i_2}\ar[dr]&&\mathcal (V/V')\otimes\mathcal O_S\\
                                                       & \mathcal L\ar[ur]&     }
$$

For (ii), it suffices to check locally at a point $p\in S/k$. The induced morphism $\mathcal L\rightarrow (V/V')\otimes\mathcal O_S$ is essentially the morphism $\mathcal W/(V'\otimes\mathcal O_S)\rightarrow (V\otimes\mathcal O_S)/(V'\otimes\mathcal O_S)$. Locally at $p$, we have
$$\xymatrix{0\rightarrow (V'\otimes\mathcal O_S)_p\ar@<.5ex>[r]^-{i_{1p}}& \mathcal W_p\ar@<.5ex>[r]^-{i_{2p}}\ar@<.5ex>[l]^-{i'_{1p}}&(V\otimes\mathcal O_S)_p\ar@<.5ex>[l]^-{i'_{2p}}}
$$
such that $i'_{1p}\circ i_{1p}=\id,\ i'_{2p}\circ i_{2p}=\id$ and $(i'_{2p}\circ i'_{1p})\circ (i_{2p}\circ i_{1p})=\id$, because $i_1,\ i_2$ and $i_2\circ i_1$ have locally free cokernels by our assumption.  Therefore,
$$ ((V\otimes\mathcal O_S)_p/(V'\otimes \mathcal O_S)_p)/(\mathcal W_p/(V'\otimes \mathcal O_S)_p)\cong ((V\otimes\mathcal O_S)_p/\mathcal W_p).
$$
Note that $((V\otimes\mathcal O_S)_p/\mathcal W_p) $ is a free module over the local ring $(\mathcal O_S)_p$, thus our lemma follows.
\hfill  $\Box$

\vspace{.5cm}

\vspace{.5cm}

\noindent The following functor $\mathfrak F_2$ is the functor of points for Description II.
\begin{Def}\label{F_2} Let $\mathfrak F_2$ be a functor:
$$\mathbf{(Sch/k)}^{\mathrm{op}}\rightarrow \mathbf{Sets},
$$
such that for any scheme $S/k$,
\begin{itemize}
  \item[(i)] $\mathfrak F_2(S)\subset\mathfrak F(S)$.
  \item[(ii)]Let $g_{i,i}:\mathcal L_i\rightarrow \mathcal (V/V_i)\otimes\mathcal O_S$ be the sub-bundle morphism in Lemma \ref{functor two definition}. When $i\neq j$, set $g_{i,j}:\mathcal L_i\rightarrow (V/V_j)\otimes \mathcal O_S$  to be the zero morphism. In this way, for any $I\subset [n]$ we get the following induced bundle morphism:
      $$ g_I:\ \ \ \mathcal V_S\bigoplus(\oplus_{i\in I}\mathcal L_i)\rightarrow \oplus_{i\in I}(V/V_i)\otimes\mathcal O_S.
      $$
      For any $I\subset[n]$, we require $\bigwedge^{r+1-d_I+|I|}g_I=0$
\end{itemize}
\end{Def}
\begin{rmk} The Functor $\mathfrak F_2$ is represented by the scheme defined in Description II. $\mathfrak F$ is represented by $\Pi_{i=1}^n \mathbb P(V/V_i)$. Let the $n$- tuple $(\mathcal W_1, \cdots, \mathcal W_n)$ be the universal $n$- tuple of sub-bundles where for each $1\leq i\leq n $, the vector bundle $\mathcal W_i$ is a sub-bundle of $V\otimes \mathcal O_{\mathbb P(V/V_i)}$ and contains $V_i\otimes \mathcal O_{\mathbb P(V/V_i)}$ as a co-rank one sub-bundle. By Lemma \ref{functor two definition}, $\mathcal W_i/V_i\otimes\mathcal O_{\mathbb P(V/V_i)}$ is $\mathcal L_i$.  Thus on $\mathbb P(V/V_i)$, the induced morphism $0\rightarrow \mathcal L_i\rightarrow (V/V_i)\otimes \mathcal O_{\mathbb P(V/V_i)}$ is the universal sub-line bundle on $\mathbb P(V/V_i)$. Plugging $\mathcal L_i$'s in Condition (ii) in Definition \ref{F_2}, and giving coordinate systems for $V$ and for $V/V_i$'s, which gives homogeneous coordinate of $\mathbb P(V/V_i)$, one gets the closed subscheme cut by the same homogeneous ideals as in Description II.
\end{rmk}

Now we arrive to the following comparison result:
\begin{Prop} The functors $\mathfrak F_1$ and $\mathfrak F_2$ are equal as sub-functors of $\mathfrak F$, thus Description I and II are equal as closed  sub-schemes.
\end{Prop}

\noindent \emph{\textbf{Proof:}}

In the whole proof, we will fix a scheme $S/k$. In order to compare the two functors, it suffices to compare on an open cover on $S$. For each point $p\in S/k$, there exits an affine open subscheme $U_p$  containing $p$, such that $(\mathcal W_1, \cdots,\mathcal W_n)|_{U_p}$ is  trivialized. For the rest of the argument, we fix such an open cover of $S$ over $k$. Henceforth we assume that the locally free sheaves $\mathcal W_i$'s are trivialised.

Since $\mathcal V_S$ and $\{ \mathcal W_i\}_{i=1}^n|_{U_p}$ are trivialized. Then we can represent $g_I$ as a matrix after choosing a basis for $\mathcal V_S|_{U_p}$ and each $(V/V_i\otimes\mathcal O_S)|_{U_p}$.

We can choose a basis for $\mathcal V_S|_{U_p}$ and $(V/V_i)\otimes\mathcal O_S|_{U_p}$, so that the quotient map $q_i:\  \mathcal V_S|_{U_p}\rightarrow (V/V_i)\otimes\mathcal O_S|_{U_p}$ can  be represented by the following matrix of the form:
$\begin{bmatrix}
A_i\\
b_i
\end{bmatrix},
$
where $A_i$ is a $(r+1-d_i)\times (r+1)$ sub-matrix and $b_i$ a $1\times (r+1)$ submatrix and furthermore the composed map $q''_i|_{U_p}: \mathcal V_S|_{U_p}\rightarrow (V/V_i)\otimes\mathcal O_S|_{U_p}\rightarrow (\mathcal V_S/\mathcal W_i)|_{U_p}$ can be represented by
$\begin{bmatrix}
A_i
\end{bmatrix}.
$

The natural sequence of maps $\mathcal W_i|_{U_p}\hookrightarrow \mathcal V_S|_{U_p}\rightarrow (V/V_i)\otimes\mathcal O_S|_{U_p}$ factors though $h_i:\ \mathcal W_i|_{U_p}/(V_i\otimes\mathcal O_S)\rightarrow (V/V_i)\otimes \mathcal O_S|_{U_p}$. Since $(\mathcal  W_i/(V_i\otimes\mathcal O_S))|_{U_p}$ is a rank one module and the basis for $(V/V_i)\otimes \mathcal O_S|_{U_p}$ has already been chosen, the matrix representation for $h_i$ is of the form :
$\begin{bmatrix}
0 \cdots 0\cdots a_i
\end{bmatrix}^T,
$
which is a $1\times (r+1-d_i)$ matrix. One can rescale the generator of $(\mathcal  W_i/(V_i\otimes\mathcal O_S))|_{U_p}$ so that the matrix for $h_i$ is of the form:
$\begin{bmatrix}
0\cdots0\cdots 1
\end{bmatrix}^T.
$

 The map $f_I|_{U_p}:\ \mathcal V_S|_{U_p}\rightarrow \oplus_{i\in I}(\mathcal V_S/\mathcal W_i)|_{U_p}$ can be represented by the matrix of the form:\\
 $f_I|_{U_P}=\begin{bmatrix}
 A_{i_1}\\
 \vdots\\
 A_{i_{|I|}}
 \end{bmatrix},
 $
 where $I=\{i_1,\cdots, i_{|I|}\}$.

 The matrix representing $g_I|_{U_p}$ is :
 $$M_I|_{U_p}=\begin{bmatrix}
 A_{i_1} &0\\
   b_{i_1} &1\\
 \vdots\\
 A_{i_k} &0&\hdots&0&\hdots\\
 b_{i_k} &0&\hdots &1&\hdots\\
 \vdots\\
 A_{i_{|I|}}&\hdots&\hdots&\hdots&0\\
 b_{i_{|I|}}&\hdots&\hdots&\hdots&1
 \end{bmatrix}
 $$

First we will show that $\mathfrak F_1$ is a subfunctor of $\mathfrak F_2$, i.e., for a scheme $S/k$, if $(\mathcal W_1,\cdots,\mathcal W_n)\in \mathfrak F_1(S)$, then it also has $(\mathcal W_1,\cdots,\mathcal W_n)\in \mathfrak F_2(S)$.

 By condition (ii) the definition of $\mathfrak F_1$, we have $\bigwedge^{r+1-d_I}f_I|_{U_p}=0$,\ i.e., every $(r+1-d_I)\times (r+1-d_I)$ minor of the matrix
 $\begin{bmatrix}
 A_{i_1}\\
 \vdots\\
 A_{i_{|I|}}
 \end{bmatrix}
 $
 is zero.  If this is true, then by elementary linear algebra, every $(r+1-d_I+|I|)\times (r+1-d_I+|I|)$ minor of the  matrix $M_I|_{U_p}$ has to be zero.

 Next, I claim that $\mathfrak F_2$ is a sub-functor of $\mathfrak F_1$.

 By condition (ii) in the definition of $\mathfrak F_2$, at the every $(r+1-d_I+|I|)\times (r+1-d_I+|I|)$ minor of the matrix $M_I|_{U_p}$ is zero, then we can choose $(r+1-d_I+|I|)\times (r+1-d_I+|I|)$ submatrix of $g_I|_{U_p}$ with the following form
 $$\begin{bmatrix}
 W_{(r+1-d_I)}&0\\
 0&E_{|I|}
 \end{bmatrix},
 $$
 where $W_{(r+1-d_I)}$ is an arbitrary $(r+1-d_I)\times (r+1-d_I)$ submatrix of $f_I|_{U_p}$ and $E_{|I|}$ is an $|I|\times |I|$ identity matrix. Then it implies that $\det W_{(r+1-d_I)}=0$, i.e., an arbitrary $(r+1-d_I)\times (r+1-d_I)$ minor of $f_I|_{U_p}$ is zero. Then we conclude that $\mathfrak F_2$ is a sub-functor of $\mathfrak F_1$.

 Therefore, we know that Description I ( i.e. $\mathcal Z_1$) and II (i.e. $\mathcal Z_2$) are equal as closed subschemes of $\Pi_{i=1}^n \mathbb P(V/V_i)$.

\hfill  $\Box$

\section{multi-degree of the closure of the image}

 In this section, we will study the multi-degree of $X(V_1,\cdots, V_n)$. Recall that
 $$M(h)=\{(m_1,\cdots, m_n)\in \mathbb Z^n_{\geq 0}: \ \ r+1-\sum_{i\in I}m_i> d_I \ \forall I\subset [n],\ \sum_{i=1}^n m_i=h\},
 $$
 and our main result of this section is :

\begin{thm}\label{multi-degree}Set $p=\max\{h; M(h)\neq\emptyset\}$. The dimension of  $X(V_1,\cdots, V_n)$ is $p$. Its multi-degree function takes value one at the integer vectors in $M(p)$ and zero otherwise.
\end{thm}

 \begin{rmk}The closure of the image of the rational map is necessarily irreducible and reduced, and the moduli functor $\mathfrak F_1$ at least contains the closure of the image (actually, $\mathfrak F_1$ is the same as $X(V_1,\cdots, V_n)$, but we will prove this later). As long as we can show that the functor $\mathfrak F_1$ is multiplicity free and the closed points of $\mathfrak F_1$ are the closed the points of $X(V_1,\cdots, V_n)$ (by Theorem \ref{closed points F_1} ),  the multi-degree of the closure of the image is the same as the multi-degree of $\mathfrak F_1$.
 \end{rmk}

Recall that multi-degree of a closed subscheme of a product of projective spaces has a classical interpretation: for $\sum_{i=1}^n m_i=p$ and if there exists a general choice of $V^i\supset V_i$ with $\codim V^i=m_i$ in $V$ such that under the closed immersion $q_c:\ \Pi_{i=1}^n\mathbb P(V^i/V_i)\hookrightarrow\Pi_{i=1}^n \mathbb P(V/V_i)$, the closed subscheme $\Pi_{i=1}^n\mathbb P(V^i/V_i)$ intersects with $X(V_1,\cdots, V_n)$ at $k$ points counted with multiplicity, then the multi-degree function takes value $k$ at $(m_1,\dots, m_n)$.

Before proving Theorem \ref{multi-degree}, we need the following supporting proposition:

\begin{Prop}\label{vector translation} Given a $(c_1,\cdots,c_n)\in \mathbb Z^n_{\geq 0}$ with $r+1-\sum_{i\in I}c_i> d_I$, for any $I\subset[n]$, there exists a general choice of $V^i\supsetneqq V_i$ where $\codim V^i=c_i$ for each $i$ , such that\  $\dim {\cap_{i\in I}V^i}=r+1-\sum_{i\in I}c_i$, for any $I\subset [n]$, and $\dim (\cap_{i=1}^n V^i)\bigcap  V_k)=\max_{I\subset [n]}\{\dim ((\cap_{i\in I}V_i)\bigcap V_k)-\sum_{i\notin I}c_i\}$ for each $1\leq k \leq n$.  In particular under this general choice of $\{V^i\}_{i=1}^n$, we have $(X(V_1,\cdots, V_n)\bigcap\Pi_{i=1}^n \mathbb P(V^i/V_i))\neq \emptyset$.
\end{Prop}

The next lemma is fundamental for the whole theory to the work out:
\begin{Lemma}\label{fundamental}Let $\{c_i\}_{i=1}^h$ be a finite sequence of nonnegative integers and $V'$ be a fixed subspace of $V$ with $\dim V=r+1$. If for any $I\subset [h]$,  we have $(r+1-\sum_{i\in I}c_i)>\dim (\cap_{i\in I}V_i)$, then there exists $V^i\supset V_i$, for each $i\in [h]$ with $\codim V^i=c_i$,  satisfying the following properties:
\begin{itemize}
\item[(i)]$\dim (\cap_{i\in [h]}V^i)=r+1-\sum_{i\in [h]}c_i$,
\item[(ii)]$\dim(\cap_{i\in [h]}V^i\ \bigcap\  V')\leq\max_{I\subset [h]}\{\dim \cap_{i\in I}V_i\ \bigcap\  V'-\sum_{i\notin I} c_i \}$.
\end{itemize}
\end{Lemma}

\noindent\emph{\textbf{Proof:}}
This proof is by induction on the cardinality of $[h]$, when $h=1$, we can find $V^1\supset V_1$, with $\codim V^1=c_1$ and $\dim (V^1\cap V')=\max\{\dim(V'\cap V_1),\dim (V')-c_1 \}$ as follows:

Consider a decomposition $V=V'\oplus V''$,such that $V_1=(V'\cap V_1)\oplus (V''\cap V_1)$. Note that
$$\dim (V'\cap V_1)\leq \dim(V_1)< r+1-c_1.
$$
Consider another decomposition $V''=(V''\cap V_1)\oplus W_{(3)}$ and $V'=(V'\cap V_1)\oplus W_{(4)}$.

When $\dim W_{(3)}\geq (r+1-c_1-\dim (V_1))$, then take a $(r+1-c_1-\dim (V_1))$ dimensional subspace of $W_{(3)}$ and denote this subspace as $W'_{(3)}$. Take $V_1\oplus W'_{(3)}$ to be $V^1$, and then one can see that $\dim (V^1\cap V')=\dim (V'\cap V_1)$.

When $\dim W_{(3)}< (r+1-c_1-\dim (V_1))$, note that $\dim W_{(3)}=r+1-\dim V'-(\dim V_1-\dim (V'\cap V_1))$, then we have $r+1-\dim V'-(\dim V_1-\dim (V'\cap V_1))<r+1-c_1-\dim (V_1)$ thus  $\dim V'-c_1> \dim (V_1\cap V')$. Take a $\dim(V')-c_1$ dimensional subspace of $V'$ containing $V_1\cap V'$ and denote this space as $W'_{(4)}$, then $\dim (W'_{(4)}\oplus V'')=(\dim (V')-c_1)+(r+1-\dim V')=r+1-c_1$. Take $W'_{(4)}\oplus V''$ to be $V^1$, and one can see that $\dim (V^1\cap V')=\dim (V')-c_1$. Thus we conclude that there exists $V^1$ with $\codim V^1=c_1$ and $\dim (V^1\cap V')=\max \{\dim(V'\cap V_1), \dim (V')-c_1\}$.

When $h\leq k$, assume the lemma is true. By the assumption, there exists a system of linear subspaces of $\{V^i\}_{i=1}^k$ satisfying the properties (i) and (ii) with respect to each of $V_{k+1},\ V'$ or $(V_{k+1}\cap V')$.

Let $\Gr(c_i, V/V_i)$  be the irreducible   moduli scheme parameterizing all codimension $c_i$ subspaces of $V$ containing $V_i$.   The subset of  $\Pi_{i=1}^k\Gr (c_i, V/V_i)$ satisfying Properties (i) and (ii)  with respect to all $V_{k+1},\ V'$ and $(V_{k+1}\cap V')$ is the intersection of three non-empty open subset of $\Pi_{i=1}^k\Gr (c_i, V/V_i)$, which is irreducible, so each open subset of this moduli space is open and dense. Note that intersection of three dense nonempty open subsets is nonempty, so there exists a system of linear subspaces $\{V^i\}_{i=1}^k$ satisfying Properties (i) and (ii) with respect to both $V_{k+1},V'$ and $V_{k+1}\cap V'$, and for the rest of the proof, we will fix such a system $\{V^i\}_{i=1}^k$.

First,  I claim that there exists $V^{k+1}\supset V_{k+1}$ with $\codim V^{k+1}=c_{k+1}$, such that $\dim \cap_{i=1}^{k+1}V^{i}=r+1-\sum_{i=1}^{k+1}c_i$. We can construct such a $V^{k+1}$ as follows:

Consider a decomposition $V=\cap_{i=1}^kV^i\oplus V^K$, such that $$V_{k+1}=((\cap_{i=1}^kV^i)\bigcap V_{k+1})\oplus (V^K\cap V_{k+1}).$$

Then $\dim V^K= (r+1)-(r+1-\sum_{i=1}^kc_i)=\sum_{i=1}^k c_i$, but $\dim V^{k+1}=r+1-c_{k+1}>\sum_{i=1}^k c_i$ (because $r+1-\sum_{i=1}^{k+1}c_i>\dim(\cap_{i=1}^{k+1}V_i)$).

$\dim(\cap_{i=1}^k V^i)=r+1-\sum_{i=1}^k c_i\geq r+1-\sum_{i=1}^{k+1}c_i$, and also notice that for any $I\subset [k]$, we have

$$\begin{array}{rcl}

\dim(\cap_{i\in I}V_i\bigcap{V_{k+1}})-\sum_{i\notin I,\ i\in [k]}c_i  &<&r+1-\sum_{I\cup\{k+1\}}c_i-\sum_{i\notin I,\ i\in [k]}c_i\\
                                                            &=&r+1-\sum_{i=1}^{k+1}c_i
\end{array}
$$

Therefore,
$$\dim(\cap_{i=1}^kV^i\ \bigcap V_{k+1})\leq\max_{I\subset[k]}\{\dim(\cap_{i\in I}V_i\bigcap{V_{k+1}})-\sum_{i\notin I}c_i\}<r+1-\sum_{i=1}^{k+1}c_i,$$
then take a $(r+1-\sum_{i=1}^{k+1}c_i)$ dimensional subspace of $\cap_{i=1}^kV^i$ containing $\cap_{i=1}^kV^i\ \bigcap V_{k+1}$, and denote this subspace as $V_{[k]}'$.
$$\begin{array}{rcl}
\dim (V^K\oplus V'_{[k]})&=&\dim (V^K)+\dim V'_{[k]}\\
                         &=&\sum_{i=1}^k c_i+(r+1-\sum_{i=1}^{k+1}c_i)\\
                         &=&r+1-c_{k+1},
\end{array}$$
then take $V^{k+1}=V^K\oplus V_{[k]}'$. Then this $V^{k+1}$ has the property that $V^{k+1}\supset V_{k+1}$ with codimension $c_k$ and $\dim (V^{k+1}\bigcap\ \cap_{i=1}^kV^i )=r+1-\sum_{i=1}^{k+1}c_i$, which is the minimal intersection dimension of $V^i$'s satisfying $\codim V^i=c_i$.

Secondly, I claim that we can construct a $V^{k+1}$ with the property $$\dim (\cap_{i=1}^{k+1}V^{i}\bigcap\ V')\leq\max_{I\subset [k+1]}\{\dim(\cap_{i\in I}V_i\ \bigcap\ V')-\sum_{i\notin I [k+1]}c_i\}.$$

Consider a decomposition $V=(\cap_{i=1}^kV^k\bigcap V')\oplus W$ such that $V_{k+1}=(\cap_{i=1}^kV^k\bigcap V'\bigcap V_{k+1})\oplus (W\cap V_{k+1})$.

When $\dim W\geq r+1-c_{k+1}-\dim (\cap_{i=1}^kV^i\bigcap V'\bigcap V_{k+1})$, take $[r+1-c_{k+1}-\dim (\cap_{i=1}^kV^i\bigcap V'\bigcap V_{k+1})]$ dimensional subspace of $W$ containing $W\cap V_{k+1}$ (one can do this because $\dim (W\cap V_{k+1})+\dim (\cap_{i=1}^kV^i\bigcap V'\bigcap V_{k+1})=\dim(V_{k+1})$ and $\dim V_{k+1}< r+1-c_{k+1}$) and denote this subspace as $W'$. Take $V^{k+1}$ to be $\spa(V_{k+1}, W')$. One can check that $\dim (\spa(V_{k+1}, W'))=r+1-c_{k+1}$, and $\dim (\spa(V_{k+1}, W')\bigcap (\cap_{i=1}^kV^i\bigcap V'))=\dim (V_{k+1}\bigcap (\cap_{i=1}^kV^i\bigcap V')) $. Then in this case there exists $V^{k+1}$, such that  $\dim (V^{k+1}\bigcap (\cap_{i=1}^k V^i\bigcap V'))\leq \max\{\dim(\cap_{i=1}^kV^i\bigcap V'\cap V_{k+1}),\ \dim(\cap_{i=1}^kV^i\ \bigcap V')-c_{k+1})\}$.

When $\dim W< r+1-c_{k+1}-\dim (\cap_{i=1}^kV^i\bigcap V'\bigcap V_{k+1})$, note that $\dim W=r+1-\dim (\cap_{i=1}^k V^i\bigcap V')$, then we have
$$r+1-\dim (\cap_{i=1}^k V^i\bigcap V')<r+1-c_{k+1}-\dim (\cap_{i=1}^kV^i\bigcap V'\bigcap V_{k+1})
,$$
thus $\dim (\cap_{i=1}^kV^i\bigcap V'\bigcap V_{k+1})<\dim (\cap_{i=1}^k V^i\bigcap V')-c_{k+1}$. Take $\dim (\cap_{i=1}^k V^i\bigcap V')-c_{k+1}$ dimensional linear subspace of $\cap_{i=1}^k V^i\bigcap V'$ containing $\cap_{i=1}^kV^i\bigcap V'\bigcap V_{k+1}$ and denote this linear space as $W''$.  Notice that
$$\begin{array}{rcl}
\dim (W''\oplus W)&=&(\dim (\cap_{i=1}^k V^i\bigcap V'))-c_{k+1}+(r+1-\dim (\cap_{i=1}^k V^i\bigcap V'))\\
                  &=&r+1-c_{k+1}
\end{array}.
$$
Also $\dim((W''\oplus W)\bigcap (\cap_{i=1}^kV^i\bigcap V'))=\dim (W''\bigcap(\cap_{i=1}^kV^i\bigcap V')))=\dim (\cap_{i=1}^k V^i\bigcap V')-c_{k+1}$. Take $V^{k+1}$ to be $W''\oplus W$, Then in this case there exists $V^{k+1}$, such that  $\dim (V^{k+1}\bigcap (\cap_{i=1}^k V^i\bigcap V'))\leq \max\{\dim(\cap_{i=1}^kV^i\bigcap V'\cap V_{k+1}),\ \dim(\cap_{i=1}^kV^i\ \bigcap V')-c_{k+1})\}$.

Combining the above two cases,  we have proved that there always exists $V^{k+1}$ whose codimension is $c_{k+1}$ such that  the intersection dimension of $V^{k+1}$ with $\cap_{i=1}^kV^i\ \bigcap\ V'$ is not bigger than
$$\max\{\dim(\cap_{i=1}^kV^i\bigcap V'\cap V_{k+1}),\ \dim(\cap_{i=1}^kV^i\ \bigcap V')-c_{k+1})\}.$$
Recall that
$$\dim(\cap_{i=1}^kV^i\bigcap V'\cap V_{k+1})\leq\max_{I\subset [k]}\{\dim(\cap_{i\in I}V_i\bigcap (V'\cap V_{k+1}))-\sum_{i\in [k]\setminus I}c_i\},$$ therefore we have:
$$
\begin{array}{rl}
 &\max\{\dim(\cap_{i=1}^kV^i\bigcap V'\cap V_{k+1}),\ \dim(\cap_{i=1}^kV^i\ \bigcap V')-c_{k+1}\}\\
\leq&\max_{I\subset[k]}\{\dim(\cap_{i\in I}V_i\bigcap (V'\cap V_{k+1}))-\sum_{i\in [k]\setminus I}c_i,(\dim \cap_{i\in I}V_i\bigcap V'-\sum_{i\in [k]\setminus I}c_i)-c_{k+1}\}\\
=&\max_{I\subset [k+1]}\{\dim(\cap_{i\in I}V_i\cap V')-\sum_{i\in [k+1]\setminus I}c_i\}
\end{array}
$$

Therefore, all of the possible $V^{k+1}$'s satisfying Case (i) form an nonempty open subset on $\Gr(c_{k+1}, V/V_{k+1})$, where $\Gr(c_{k+1}, V/V_{k+1})$ is the irreducible moduli scheme parameterizing all of the codimension $c_{k+1}$ linear subspaces in $V$ containing $V_{k+1}$. All of the possible $V^{k+1}$ satisfying Case (ii)  also form  an nonempty open subset on $\Gr(c_{k+1}, V/V_{k+1})$. So the all of the possible $V^{k+1}$ satisfying both Cases (i) and (ii) form the intersection of two nonempty open subsets of the $\Gr(c_{k+1}, V/V_{k+1})$, which is still open and dense because $\Gr(c_{k+1}, V/V_{k+1})$ is irreducible, thus nonempty, so there exists one $V^{k+1}$ satisfying both Cases (i) and (ii).

Thus we have found $\{V^i\}_{i=1}^{k+1}$ satisfying Properties (i) and (ii). Then the lemma follows.

\hfill $\Box$

\begin{rmk}\label{why mini} A priori, given $V^i\supset V_i$ with $\codim V^i=c_i$ for each $i$, we always have $\dim (\cap_{i=1}^n V^i\bigcap V')\geq \dim (\cap_{i\in I}V_i\bigcap V')-\sum_{i\notin I}c_i$ for any $I\subset [n]$, therefore $\max_{I\subset [h]}\{\dim \cap_{i\in I}V_i\ \bigcap\  V'-\sum_{i\notin I} c_i \}$ in  the Lemma \ref{fundamental} \emph{is} the actual minimal possible intersection dimension of $\bigcap_{i\in [n]}V^i$ with $V'$
\end{rmk}

Now we are in a good position to prove Proposition \ref{vector translation}:

\noindent \emph{\textbf{Proof of Proposition \ref{vector translation}}}:

Take $V'=V_k,\ 1\leq k\leq n $,  then by Lemma \ref{fundamental}, we can find a sequence $\{V^i\supset V_i\}_{i=1}^n$, such that $\dim (\cap_{i\in I}V^i)=r+1-\sum_{i\in I} c_i$ for any $I\subset [n]$, also $\dim ((\cap_{i=1}^nV^i)\bigcap V_k)\leq\max_{I\subset [n]}\{\dim((\cap_{i\in I}V_i)\bigcap V_k)-\sum_{i\in[n]\setminus I}c_i\}$. By Remark \ref{why mini}, $\max_{I\subset [h]}\{\dim \cap_{i\in I}V_i\ \bigcap\  V_k-\sum_{i\notin I} c_i \}$ in  the Lemma \ref{fundamental} \emph{is} the actual minimal possible intersection dimension of $\bigcap_{i\in [n]}V^i$ with $V_k$.

Note that $\{V^i\}_{i=1}^n$ with $\codim V^i=c_i$ for each $i$ corresponds to a point in the following irreducible moduli scheme
 $\Pi_{i=1}^n \Gr (c_i, V/V_i)$ where $\Gr (c_i, V/V_i)$ is the moduli scheme parameterizing all of codimension $c_i$ linear subspaces containing $V_i$. The condition $\dim \cap_{i\in I}V^i=r+1-\sum_{i\in I} c_i$ for any $I\subset [n]$ is a nonempty open condition  and also $\dim (\cap_{i=1}^n V^i =1)\bigcap V_k=\max_{I\subset [h]}\{\dim \cap_{i\in I}V_i\ \bigcap\  V_k-\sum_{i\notin I} c_i \}$ is also a nonempty open condition by the above discussion.  Therefore for each $k\in [n]$, satisfying $\dim\cap_{i\in I} V^i =r+1-\sum_{i\in I}c_i$ and $\dim (\cap_{i=1}^n V^i =1)\bigcap V_k=\max_{I\subset [h]}\{\dim \cap_{i\in I}V_i\ \bigcap\  V_k-\sum_{i\notin I} c_i \} $ is the intersection of two nonempty open dense subsets of the moduli scheme, which is still a nonempty open condition in $\Pi_{i=1}^n \Gr (c_i, V/V_i)$, and we denote this nonempty open set of $\Pi_{i=1}^n \Gr (c_i, V/V_i)$ as $U_k$.

Then we take $\cap_{i=1}^n U_k$ is still nonempty and open, thus we can take a point in this nonempty set, which corresponds to $\{V^i\}_{i=1}^n$ with $\dim(\cap_{i\in I} V^i)=r+1-\sum_{i\in I}c_i$ and $\dim((\cap_{i=1}^n V^i)\bigcap (V_k))=\max_{I\subset [h]}\{\dim \cap_{i\in I}V_i\ \bigcap\  V_k-\sum_{i\notin I} c_i \}$.

Since $d_I< r+1-\sum_{i\in I}c_i$, for any $I\subset [n]$, therefore
$$\begin{array}{rcl}
 \dim((\cap_{i\in [n]}V^i)\bigcap V_k)&=& \max_{I\subset [n]}\{\dim ((\cap_{i\in I}V_i)\bigcap V_k)-\sum_{i\notin I}c_i\}\\
                                      &<& r+1-\sum_{i=1}^n c_i\\
                                      &=&\dim\bigcap_{i=1}^n V^i   .
 \end{array}
$$

Then for each $i$, there exists a vector $v_i\in \cap_{j\in [n]}V^j$ such that $v_i\notin V_i$, then there exists a vector $v_{[n]}\in \cap_{i\in [n]}V^i$ such that $v_{[n]}\notin V_i$ for each $i$.  Therefore, the intersection of $\Pi_{i=1}^n \mathbb P(V^i/V_i)$ with $X(V_1,\cdots, V_n)$ contains a point $(\spa(V_1, v_{[n]}),\cdots, \spa(V_n, v_{[n]}))$, which implies for such a general choice of $\{V^i\}_{i=1}^n$, the intersection is not empty.

\hfill $\Box$

\vspace{0.3cm}

The following lemma is used to detect at which integer vector the multi-degree function takes zero.
\begin{Lemma}\label{if and only if} Given an integer vector $(c_1,\cdots, c_n)$  such that there exists $I\subset [n]$ with $r+1-\sum_{i\in I}c_i\leq \dim(\cap_{i\in I}V_i)$, then there exists $\Pi_{i=1}^n \mathbb P^{r+1-d_i-c_i}\hookrightarrow \Pi_{i=1}\mathbb P(V/V_i)$,such that it has empty intersection with $X(V_1,\cdots, V_n)$..
\end{Lemma}

\emph{\textbf{Proof:}}

Rational points of $\Pi_{i=1}^n \mathbb P^{r+1-d_i-c_i}\hookrightarrow \Pi_{i=1}\mathbb P(V/V_i)$ correspond to $(V^1,\cdots, V^n)$ with $V^i\supsetneqq V_i$ and $\codim V^i=c_i$. If we want to show there exists $\Pi_{i=1}^n \mathbb P^{r+1-d_i-c_i}\hookrightarrow \Pi_{i=1}\mathbb P(V/V_i)$ with empty intersection with the closure of the rational map $\mathbb P(V)\dashrightarrow \mathbb P(V/V_i)$, all we need to do is to construct a system of $\{V^i\supset V_i\}_{i=1}^n$ such that it does not contain $(W_1,\dots, W_n)$ with $\forall I\subset [n], \cap_{i\in I}W_i\supsetneqq \cap_{i\in I}V_i$.

First observe that since there exists $I\subset [n]$ with $r+1-\sum_{i\in I}c_i\leq \dim(\cap_{i\in I}V_i)$, then there exists $I'\subset I$, such that for any $ I''\subsetneqq I',\ r+1-\sum_{i\in I''}c_i> \dim(\cap_{i\in I''}V_i)$ but also $r+1-\sum_{i\in I'}c_i\leq\dim(\cap_{i\in I'}V_i)$.

When $|I'|=1$, we have $r+1-c_{I'}\leq d_{I'}$. In this case, any linear subspace of $V$ strictly containing $V_{I'}$ has codimension less than $c_{I'}$. Thus any $\Pi_{i=1}^n \mathbb P^{r+1-d_i-c_i}\hookrightarrow \Pi_{i=1}\mathbb P(V/V_i)$ has empty intersection with $X(V_1,\cdots, V_n)$.

When $|I'|>1$, we will  construct $\{V^i\}_{i\in [n]}$ with $\codim V^i=c_i,\forall i\in [n]$ such that there exists $I'\subset[n]$ with $\dim(\cap_{i\in I'}V^i)=\dim(\cap_{i\in I'}V_i)$.
Let $i''\in I'$, consider $I'\setminus i''$, then by Lemma \ref{fundamental},  there exists $\{V^i\supset V_i\}_{i\in I'\setminus i''}$, with $\dim(\cap_{i\in I'\setminus i''}V^i)=r+1-\sum_{i\in I'\setminus i''}c_i$ and also $$\dim (\cap_{i\in I'\setminus i''}V^i)\bigcap V_{i''}\leq\max_{\mathcal K\subset I''}\{\dim(\cap_{\mathcal K\subset I''}V_i\ \bigcap V_{i''})-\sum_{i\in I''\setminus\mathcal K}c_i\},$$ where $I''$ is $I'\setminus i''$.

I claim  that $\dim(\cap_{i\in I''}V^i)\bigcap V_{i''}=\dim(\cap_{i\in I'}V_i)$.

First let's prove that $\max_{\mathcal K\subset I''}\{\dim(\cap_{\mathcal K\subset I''}V_i\ \bigcap V_{i''})-\sum_{i\in I''\setminus\mathcal K}c_i\}=\dim(\cap_{i\in I'}V_i)$.

When $\mathcal K\subsetneqq I''$, we have $\mathcal K\cup \{i''\}\subsetneqq I'$, then
$$\begin{array}{rcl}
\dim(\cap_{i\in\mathcal K\subset I''}V_i\ \bigcap V_{i''})-\sum_{i\in I''\setminus\mathcal K}c_i&\leq& (r+1-\sum_{\mathcal K\cup\{i''\}}c_i)-\sum_{i\in I''\setminus\mathcal K}c_i\\
                       &=&r+1-\sum_{i\in I'}c_i\\
                       &\leq&\dim(\cap_{i\in I'}V_i)
\end{array}
$$

When $\mathcal K=I''$, we have
 $\dim(\cap_{\mathcal K\subset I''}V_i\ \bigcap V_{i''})-\sum_{i\in I''\setminus\mathcal K}c_i=\dim(\cap_{i\in I'}V_i)$. Thus it follows that $\max_{\mathcal K\subset I''}\{\dim(\cap_{\mathcal K\subset I''}V_i\ \bigcap V_{i''})-\sum_{i\in I''\setminus\mathcal K}c_i\}=\dim(\cap_{i\in I'}V_i)$. Therefore $\dim(\cap_{i\in I''}V^i)\bigcap V_{i''}\leq \dim(\cap_{i\in I'}V_i)$

 But since $V^i\supset V_i$, we have $\dim(\cap_{i\in I''}V^i)\bigcap V_{i''}\geq \dim(\cap_{i\in I'}V_i)$. Then we have $\dim(\cap_{i\in I''}V^i)\bigcap V_{i''}=\dim(\cap_{i\in I'}V_i)$, thus the claim follows.

 Now consider a decomposition $V=\cap_{i\in I''}V^i\oplus V^D$ such that $V_{i''}=((\cap_{i\in I''}V^i)\bigcap V_{i''})\oplus (V^D\cap V_{i''})$.  Note that $\dim (V^D)=(r+1)-(r+1-\sum_{i\in I''}c_i)=\sum_{i\in I''}c_i\geq r+1-c_{i''}-\dim \cap_{i\in I''}V^i\bigcap V_{i''}$, so let's take $(r+1-c_{i''}-\dim \cap_{i\in I''}V^i\bigcap V_{i''})$-dimensional subspace $V^d$ of $V^D$ such that $V^d\supset V_{i''}\cap V^D$ and take $V^{i''}=\spa(V_{i''}, V^d)$. (note that $\dim (V_{i''}\cap V^D)+\dim (\cap_{i\in I''}V^i\bigcap V_{i''})=\dim V_{i''}<r+1-c_{i''}$).

 \emph{Therefore, we have proved that there exists $\{V^i\supset V_i\}_{i\in I'}$ with $\codim V^i=c_i$ and $\dim(\cap_{i\in I'}V^i)=\dim\cap_{i\in I'}V_i$.}

Also, points of the closure of the rational map $\mathbb P(V)\dashrightarrow \mathbb P(V/V_i)$ correspond to $(W_1,\cdots,W_n)$ with $W_i$ a linear subspace of $V$ containing $V_i$ as a codimension one subspace, such that $\forall \mathcal K\subset[n],\dim (\cap_{i\in \mathcal K}W_i)>\dim(\cap_{i\in \mathcal K}V_i)$.  Then the $\{V^i\}_{i=1}^n$ satisfying for a $I'\subset [n],\ \dim (\cap_{i\in I'}V^i)=\dim \cap_{i\in I}(V_i)$, does not contain any $(W_1,\cdots,W_n)$ because  $\dim(\cap_{i\in I'}W_i)>\dim (\cap_{i\in I'}V_i)$. Thus the lemma follows.

\hfill $\Box$

\begin{Prop}\label{complete multi-degree} The dimension of the image of the rational map:
$$  \mathbb P(V)\dashrightarrow \Pi_{i=1}^n \mathbb P(V/V_i),\ \ \bigcap_{i=1}^n V_i=\{0\}
$$
is the largest $p$ such that $M(p)\neq \emptyset$. And its multi-degree function takes a nonzero value at the integer vectors in $M(p)$ and zero otherwise.
\end{Prop}

\noindent \emph{\textbf{Proof:}}
In this proof, $t_i=\dim (V/V_i)-1$ for each $1\leq i\leq n$, so for the rest of the proof, we identify $\mathbb P(V/V_i)$ with $\mathbb P^{t_i}$. Given an ample divisor $\sum_{i=1}^n c_i \mathbb P^{t_1}\times \cdots \times\mathbb P^{t_i-1}\times\cdots \times \mathbb P^{t_n}$ in $\Pi_{i=1}^n\mathbb P^{t_i}\hookrightarrow \mathbb P^{\ell}$, where $c_i\gg 0$ for each $1\leq i\leq n$, we get a closed immersion:
$$f_c:\ \ \ \Pi_{i=1}^n\mathbb P^{t_i}\hookrightarrow \mathbb P^{\ell}
$$
Let $D(p)=\{(u_1,\cdots, u_n)\in \mathbb Z^n_{\geq 0}:\ \sum_{i=1}^n u_i= p\}$.  I claim that in Chow ring of $\Pi_{i=1}^n \mathbb P^{m_i}$, the class $(\sum_{i=1}^n c_i \mathbb P^{t_1}\times \cdots \times\mathbb P^{t_i-1}\times\cdots \times \mathbb P^{t_n})^k$ where $u_i\leq t_i$ for each $1\leq i\leq n$, can be represented by the form $\sum_{(u_1,\cdots, u_n)\in D(k)}c_{(u_1,\cdots, u_n)}\Pi_{i=1}^n \mathbb P^{t_i-u_i}$ where all $c_{(u_1,\cdots, u_n)}$ are positive.
One can prove this by induction. Assume that when $h=k$ it is correct. Let $e_j=(0,\cdots,1,\cdots,0)$ whose $j$-th component is one and other components are zero. When $h=k+1$,
$$ \begin{array}{rcl}
&&(\sum_{i=1}^n c_i \mathbb P^{t_1}\times \cdots \mathbb P^{t_i-1}\cdots \mathbb P^{t_n})^{k+1}\\
&=&(\sum_{i=1}^n c_i \mathbb P^{t_1}\times \cdots \mathbb P^{t_i-1}\cdots \mathbb P^{t_n})^k\cdot (\sum_{i=1}^n c_i \mathbb P^{t_1}\times \cdots \times\mathbb P^{t_i-1}\times\cdots\times \mathbb P^{t_n})\\
 &=&(\sum_{(u_1,\cdots, u_n)\in D(k)}c_{(u_1,\cdots, u_n)}\Pi_{i=1}^n \mathbb P^{t_i-u_i})\cdot (\sum_{i=1}^n c_i \mathbb P^{t_1}\times \cdots \times\mathbb P^{t_i-1}\times\cdots \times \mathbb P^{t_n})\\
 &=&(\sum_{(u_1,\cdots, u_n)\in D(k)}c_{(u_1,\cdots, u_n)}\Pi_{i=1}^n \mathbb P^{t_i-u_i})\cdot (\sum_{i=1}^n c_i \mathbb P^{t_1}\times \cdots \times\mathbb P^{t_i-1}\times\cdots\times \mathbb P^{t_n})\\
 &=&(\sum_{(u_1,\cdots, u_n)\in D(k)}c_{(u_1,\cdots, u_n)}\sum_{j=1}^n(c_j\Pi_{i\in [n]\setminus j} \mathbb P^{t_i-u_i}\times \mathbb P^{t_j-u_j-1}))\\
 &=&\sum_{(u'_1,\cdots, u'_n)\in D(k+1)})c'_{(u'_1,\cdots, u'_n)}\Pi_{i=1}^n \mathbb P^{t_i-u'_i},

\end{array}
$$

where all $c'_{(u'_1,\cdots, u'_n)}=\sum_{(u_1,\cdots, u_n)\in \mathcal M(u'_1, \cdots, u'_n)}c_{(u_1,\cdots, u_n)}\cdot c_{n(u_1,\cdots, u_n)}$ with
$$\mathcal M(u'_1,\cdots, u'_n)=\{(u_1,\cdots, u_n)\in D(k);\ \exists n(u_1,\cdots,u_n)\in \mathbb Z_{\geq 0}\  \mathrm{with}\ (u_1,\cdots, u_n)+e_{n(u_1,\cdots,u_n)}=(u'_1,\cdots,u'_n)\},$$
which clearly implies that $c'_{(u'_1,\cdots, u'_n)}$ is positive for any $(u'_1,\cdots,u'_n)\in D(k+1)$.

Note that
$$(\sum_{i=1}^n c_i \mathbb P^{t_1}\times \cdots \mathbb P^{t_i-1}\cdots \times \mathbb P^{t_n})^k={f^*_c}\mathbb P^{\ell-k},
$$
and then we have:
$$\begin{array}{rcl}
X(V_1,\cdots, V_n)\cdot {f^*_c}\mathbb P^{\ell-k}&=&X(V_1,\cdots, V_n)\cdot \sum_{(u_1,\cdots, u_n)\in D(p)}c_{(u_1,\cdots, u_n)}\Pi_{i=1}^n \mathbb P^{m_i-u_i}.\\
\end{array}
$$
Since $p$ is the largest integer with $M(p)=\{(m_1,\cdots, m_n):\ r+1-\sum_{i\in I}>d_I,\forall I\subset[n],\ \sum_{i=1}^n m_i=p\}\neq \emptyset$, we can take an element $(m_1,\cdots, m_n)\in M(p)$. By Proposition \ref{vector translation}, for a general choice of $V^i\supset V_i$ satisfying the conditions in Proposition \ref{vector translation},  $\Pi_{i=1}^n\mathbb P(V^i/V_i)\hookrightarrow \Pi_{i=1}^n \mathbb P(V/V_i)$ always has nonempty intersection with $X(V_1,\cdots, V_n)$.  For $(u_1,\cdots,u_n)\in D(p)\setminus M(p)$, by Lemma \ref{if and only if}, for a generic choice of $V^i \subset V$ with $\codim V^i=j_i$, the closed subvariety $\Pi_{i=1}^n\mathbb P(V^i/V_i)\hookrightarrow \Pi_{i=1}^n \mathbb P(V/V_i)$  has empty intersection with $X(V_1,\cdots, V_n)$. Then we have:.
$$\begin{array}{rcl}
X(V_1,\cdots, V_n)\cdot {f^*_c}\mathbb P^{\ell-p}&=&X(V_1,\cdots, V_n)\cdot \sum_{(u_1,\cdots, u_n)\in D(p)}c_{(u_1,\cdots, u_n)}\Pi_{i=1}^n \mathbb P^{t_i-u_i}\\
               &=& X(V_1,\cdots, V_n)\cdot \sum_{(u_1,\cdots, u_n)\in M(p)}c_{(u_1,\cdots, u_n)}\Pi_{i=1}^n \mathbb P^{t_i-u_i}\\
               &=&\sum_{(u_1,\cdots, u_n)\in M(p)}c_{(u_1,\cdots, u_n)}(\Pi_{i=1}^n \mathbb P^{t_i-u_i}\cdot X(V_1,\cdots, V_n)).
\end{array}
$$

Note that
$$\begin{array}{rcl}
 & &{f_c}_*(\sum_{(u_1,\cdots, u_n)\in M(p)}c_{(u_1,\cdots, u_n)}(\Pi_{i=1}^n \mathbb P^{m_i-u_i}\cdot X(V_1,\cdots, V_n)))\\
 &=&(\sum_{(u_1,\cdots, u_n)\in M(p)}c_{(u_1,\cdots, u_n)}{f_c}_*(\Pi_{i=1}^n \mathbb P^{t_i-u_i}\cdot X(V_1,\cdots, V_n)).
\end{array}
$$
 By the projection formula, the intersection class ${f_c}_*(X(V_1,\cdots, V_n))\cdot \mathbb P^{\ell-p}$ is also \\
 $(\sum_{(u_1,\cdots, u_n)\in M(p)}c_{(u_1,\cdots, u_n)}{f_c}_*(\Pi_{i=1}^n \mathbb P^{t_i-u_i}\cdot X(V_1,\cdots, V_n))$.  The intersection product $\Pi_{i=1}^n \mathbb P^{t_i-u_i}\cdot X(V_1,\cdots, V_n)$ corresponds to a nontrivial closed subvariety of $\Pi_{i=1}^n \mathbb P (V/V_i)$, which is in turn a closed subvariety of $\mathbb P^{\ell}$ under the closed immersion $f_c$. Hence  ${f_c}_*(\Pi_{i=1}^n \mathbb P^{t_i-u_i}\cdot X(V_1,\cdots, V_n))$ corresponds to a closed subvariety of $\mathbb P^{\ell}$ and $c_{(u_1,\cdots,u_n)}$ is a positive integer for each $(u_1,\cdots, u_n)\in M(p)$. Therefore, $(\sum_{(u_1,\cdots, u_n)\in M(p)}c_{(u_1,\cdots, u_n)}{f_c}_*(\Pi_{i=1}^n \mathbb P^{m_i-u_i}\cdot X(V_1,\cdots, V_n))$ is a nontrivial class in the Chow ring of $\mathbb P^{\ell}$, which implies that the codimension $p$ closed subvariety $\mathbb P^{\ell-p}$ always has nontrivial intersection with ${f_c}_*(X(V_1,\cdots, V_n))$. Hence we have $r'\geq p$, where $r'$ is the dimension of $X(V_1,\cdots, V_n)$.

Note that $M(p+1)=\emptyset$, so for any $(j_1,\cdots, j_n)$ with $\sum_{i=1}^nj_i=p+1$, there exists $I\subset [n]$, with $r+1-\sum_{i\in I}j_i\leq d_I$, then by Lemma \ref{if and only if}, for a generic choice of $V^i \subset V$ with $\codim V^i=j_i$, the closed subvariety $\Pi_{i=1}^n\mathbb P(V^i/V_i)\hookrightarrow \Pi_{i=1}^n \mathbb P(V/V_i)$  has empty intersection with $X(V_1,\cdots, V_n)$. Therefore we have:
$$\begin{array}{rcl}
X(V_1,\cdots, V_n)\cdot {f^*_c}\mathbb P^{\ell-(p+1)}&=&X(V_1,\cdots, V_n)\cdot \sum_{(u_1,\cdots, u_n)\in D(p)}c_{(u_1,\cdots, u_n)}\Pi_{i=1}^n \mathbb P^{m_i-u_i}\\
               &=&0.
\end{array}
$$
 Therefore, $({f_c}_*X(V_1,\cdots, V_n))\cdot \mathbb P^{\ell-k}$ is also trivial by the projection formula, which implies the codimension $p$ closed subvariety $\mathbb P^{\ell-(p+1)}$ always has empty intersection with ${f_c}_*(X(V_1,\cdots, V_n))$. Hence we have  $r'< (p+1)$. Combining with the fact $r'\geq p$, we conclude
$$ r'=\dim(X(V_1,\cdots, V_n))=p
$$

In summary,  we have proven the domain of the multi-degree function is $D(p)=\{(m_1,\cdots, m_n)\in \mathbb Z^n_{\geq 0}: \sum_{i=1}^n m_i=p\}$.
By Lemma \ref{if and only if}, for any $(h_1,\cdots, h_n)\in D(p)\setminus M(p)$, there exists $V^i\supset V_i$ with $\codim V^i=h_i$,such that $\Pi_{i=1}^n\mathbb P(V^i/V_i)\hookrightarrow \Pi_{i=1}^n \mathbb P(V/V_i)$ has empty intersection with $X(V_1,\cdots, V_n)$, therefore multi-degree function takes zero at $D(p)\setminus M(p)$.

As for any $(m_1,\cdots, m_n)\in M(p)$, for any choice of $V^i\supset V_i$ with $\codim V^i=m_i$, the closed subvariety $\Pi_{i=1}^n\mathbb P(V^i/V_i)\hookrightarrow \Pi_{i=1}^n \mathbb P(V/V_i)$ always has nonempty intersection with $X(V_1,\cdots, V_n)$, so we conclude that multi-degree function will take a nonzero value in $M(p)$.

\hfill   $\Box$

Finally, we need to answer if multi-degree function takes a nonzero value at $(m_1,\cdots,m_n)$, then what is the nonzero value?
Let's first recall an elementary lemma about vector bundles on $X/k$ where $X/k$ is a scheme over $k$.

\begin{Lemma} \label{vector bundle} Let $X/k$ be a scheme of finite type over $k$, and let
$$\begin{xymatrix}{0\ar[r]&\mathcal F_1\ar[r]^-f& \mathcal F_2}
\end{xymatrix}
$$
be a short exact sequence of vector bundles on $X/k$. If $\coker f$ is also locally free and $\rk \mathcal F_1=\rk \mathcal F_2$, then $f$ is an isomorphism.
\end{Lemma}
\emph{\textbf{Proof:}}

It suffices to check this on the level of stalks. Let $p\in X/k$ be a point of $X/k$, notice that affine locally that exact sequence holds and localization is an exact functor, thus we have the following stalk locally at $p$:
$$\begin{xymatrix}{0\ar[r]&(\mathcal {F}_1)_p\ar[r]^{f_p}& \mathcal ({F}_2)_p\ar[r]&(\coker f)_p\ar[r]&0}
\end{xymatrix}
$$

Also we have $(\coker f)_p$ is a free module over a local ring, thus is projective. Therefore we have a splitting,
$  \begin{xymatrix}{(\mathcal {F}_1)_p\ar@<.5ex>[r]^{f_p}& \mathcal ({F}_2)_p \ar@<.5ex>[l]^{f'_p}}
\end{xymatrix}
$ with $f'_p\circ f_p=\id_{(\mathcal {F}_1)_p}$, thus $\xymatrix{\mathcal {F}_1|_p\ar[r]^{f|_{p}}& \mathcal {F}_2|_p}$ is also an injection. Let's consider the following right exact sequence:
$$\xymatrix{\mathcal {F}_1|_p\ar[r]^{f|_p}& \mathcal {F}_2|_p\ar[r]&\coker f|_p\ar[r]&0}$$
where $f|_p$ is an injection. But we have $\rk \mathcal F_1=\rk\mathcal F_2$, which implies  $\dim(\mathcal {F}_1|_p)
 =\dim(\mathcal {F}_2|_p)$ as $k_o$ vector spaces, thus $\coker f|_p=\{0\}$. By Nakayama's lemma,$(\coker f)_{p}=0$ for any $p\in X/k$, so $\coker f=0$. Therefore we conclude that $f:\ \mathcal F_1\rightarrow \mathcal F_2$ is an isomorphism.

\hfill $\Box$

\vspace{.3cm}

The following functor we are going to define represents the intersection of $X(V_1,\cdots, V_n)$ with $\Pi_{i=1}^n \mathbb P(V^i/V_i)$ in $\Pi_{i=1}^n \mathbb P(V/V_i)$, where $\{V^i\}_{i=1}^n$ were defined in Proposition \ref{vector translation}.
\begin{Def}\label{intersection functor}Let $\mathfrak H$ be a functor:
$$\mathbf{(Sck/k)}^{\mathrm{op}}\rightarrow \mathbf{Sets},
$$
and for each $S/k$, we associate the tuples of vector sub-bundles $(\mathcal W_1,\cdots,\mathcal W_n)$ of $V\otimes\mathcal O_S$ and   $\rk(\mathcal W_i)=\dim V_i+1$ for each $i\in [n]$ with the following properties:
\begin{itemize}
  \item[(a)] for any $I\subset [n]$ and $g_{I}:\ (V/\cap_{i\in I}V_i)\otimes \mathcal O_S\rightarrow \bigoplus_{i\in I}(V\otimes\mathcal O_S)/\mathcal W_i$, we have $\bigwedge^{r+1-d_I}g_{I}=0$
  \item[(b)] $\mathcal W_i$ is a sub-bundle of $V^i\otimes\mathcal O_S$ for each $i\in [n]$.
  \item[(c)]for each $i\in [n]$, we have $\xymatrix{0\rightarrow V_i\otimes\mathcal O_S\ar[r]^-{f_i}&\mathcal W_i\ar[r]^-{g_i}& \mathcal V}$ where $f_i,\ g_i$ and $g_i\circ f_i$ all have locally free cokernels.
\end{itemize}
\end{Def}

\begin{Prop}\label{answer multiplicity}The functor $\mathfrak H$ is a one point functor and thus, the multi-degree is one.
\end{Prop}

\noindent \emph{\textbf{Proof:}}

Recall at the end of the Proof of Proposition \ref{vector translation}, we have $v_{[n]}\in\bigcap_{i=1}^nV^i\setminus \bigcup_{i=1}^nV_i$.
It suffices to prove that $\mathcal W_i=\spa(V_i, v_{[n]})\otimes \mathcal O_S$ for each $0\leq i\leq n$ on any $S/k$. For the rest of the proof, let's fix a scheme $S/k$.

For any $i\in [n]$, define $J_i=\{I\subset [n]:\ i\in I\}$.  For $(m_1,\cdots, m_n)\in M(p)$, set $k_i=\min\{(r+1-\sum_{i\in I}m_i-d_I),\ \forall I\in J_i\}$. A priori, $k_i\geq 1$ and for each $i\in [n]$, we have $m_i\leq r-d_i$.

I claim that $k_i=1$.   If not, consider $(m'_1,\cdots, m'_i,\cdots, m'_n)$ where
$$m'_p=\left\{
\begin{array}{cl}m_p,& p\neq i\\
                 m_p+1,&p=i
\end{array}
\right.
$$
$(m'_1,\cdots, m'_n)$ also satisfies for any $I\subset [n],\ r+1-\sum_{i\in I}m'_i>d_I$, but $\sum_{i=1}^n m'_i=p+1$, which contradicts that $p$ is the maximal number such that $M(p)\neq \emptyset$.

Therefore, $k_i=1$, then there exits $I_i\subset [n]$ where $i\in I_i$ such that $r+1-\sum_{k\in I_i}m_k=d_{I_i}+1$.

Since $\cap_{k\in I_i}V_k\otimes \mathcal O_S$ is a sub-bundle of $\mathcal W_k$ for each $k\in I_i$, and each $\mathcal W_k$ is a sub-bundle of $V^k\otimes\mathcal O_S$ then for each $k\in I_i$, the bundle morphism $\cap_{k\in I_i}V_k\otimes\mathcal O_S\rightarrow V^k\otimes \mathcal O_S$ factors through $\mathcal W_k$. Thus the natural morphism
$$ (V/\cap_{k\in I_i}V_k)\otimes\mathcal O_s\rightarrow \bigoplus_{k\in I_i}(V/V^k)\otimes \mathcal O_S
$$
factors through $\bigoplus_{k\in I_i}((V\otimes\mathcal O_S)/(\mathcal W_k))$.
Therefore, we have the following natural morphism:
$$\xymatrix{(V/\cap_{k\in I_i}V_k)\otimes\mathcal O_s\ar[r]^-{g_{I_i}}&\bigoplus_{k\in I_i}((V\otimes\mathcal O_S)/(\mathcal W_k))\ar[r]^-{g_{w_{I_i}}}&\bigoplus_{i\in I_i}(V/V^k)\otimes \mathcal O_S},
$$

Note that $\bigwedge^{r-d_{I_i}}g_{I_i}$ is a unit ideal sheaf ($\bigwedge^{r-d_{I_i}}g_{w_{I_i}\circ g_{I_i}}$ is a unit ideal sheaf) and by the definition of $\mathfrak H$,  we have $\bigwedge^{r+1-d_{I_i}}g_{I_i}=0$. Then  by Prop. 20.8 \cite{eisen} , $\ker g_{I_i}$ is locally free. Hence for each $p\in S/k$, we have
$$(\ker g_{I_i})|_p\cong \ker (g_{I_i}|_p)=(\cap_{k\in I_i} \mathcal W_k|_p)/(\cap_{k\in I_i}V_k\otimes \mathcal O_S|_p). $$

Notice that the natural morphism

and $\ker (g_{w_{I_i}}\circ g_{I_i})=\left((\cap_{i\in I}V^i)/(\cap_{k\in I_i}V_k)\right)\otimes \mathcal O_S$, then  there is a natural sheaf injection
$$h_{I_i}:\ \ \ker g_{I_i}\rightarrow \ker (g_{w_{I_i}}\circ g_{I_i})=\left((\cap_{k\in I_i}V^k)/(\cap_{k\in I_i}V_k)\right)\otimes \mathcal O_S.$$

$h_{I_i}$ is an injection of line bundles. Next I will prove that $h_{I_i}$ is an isomorphism. It suffices to show that $\coker h_{I_i}=0$. By a geometric Nakayama argument, it suffices to show that $(\coker h_{I_i})|_p=0$ for each $p\in S/k$. Since there is a natural surjection (actually, isomorphism):
$$(\coker h_{I_i}|_p)\rightarrow (\coker h_{I_i})|_p\rightarrow 0, $$
it suffices to show $\coker h_{I_i}|_p=0$.

Fiber-wise at $p\in S/k$,  for $g_{I_i}|_p:\ V/\cap_{k\in I_i}V_k\otimes \mathcal O_S|_p\rightarrow \left(\bigoplus_{k\in I_i}\mathcal V_S/\mathcal W_k\right)|_p $, we have $\bigwedge^{r+1-d_{I_i}}=0$, this implies $\cap_{k\in I_i}V_k\otimes \mathcal O_S|_p\subsetneqq \cap_{k\in I_i}\mathcal W_k|_p$, which implies that $d_{I_i}< \dim(\cap_{k\in I_i}\mathcal W_k|_p)$. A priori, $\cap_{k\in I_i}\mathcal W_k|_p\subset \cap_{k\in I_i}V^k\otimes\mathcal O_S|_p$, which implies $\dim(\cap_{k\in I_i}\mathcal W_k|_p)\leq \dim (\cap_{k\in I_i}V^k\otimes\mathcal O_S|_p)=d_I+1$. Therefore,  $\dim(\cap_{k\in I_i}\mathcal W_k|_p)=\dim (\cap_{k\in I_i}V^k\otimes\mathcal O_S|_p)$, so we conclude that $\cap_{k\in I_i}\mathcal W_k|_p=\cap_{k\in I_i}V^i\otimes\mathcal O_S|_p$, and therefore
$$(\cap_{k\in I_i} \mathcal W_k|_p)/(\cap_{k\in I_i}V_k\otimes \mathcal O_S|_p)=\left((\cap_{k\in I_i}V^k)/(\cap_{k\in I_i}V_k)\right)\otimes \mathcal O_S|_p.
$$
Hence we have that $\coker h_{I_i}|_p=0$, thus $(\coker h_{I_i})|_p=0$ for any $p\in S/k$, therefore we conclude that $h_{I_i}$ is an isomorphism of line bundles. Equivalently, we have:
$$\ker g_{I_i}=\ker (g_{w_{I_i}}\circ g_{I_i})=\left((\cap_{k\in I_i}V^k)/(\cap_{k\in I_i}V_k)\right)\otimes \mathcal O_S
$$

 Since $\ker g_{I_i}=\left((\cap_{k\in I_i}V^k)/(\cap_{k\in I_i}V_k)\right)\otimes \mathcal O_S$ and $\mathcal W_k$ contains $(\cap_{k\in I_i}V^k)\otimes\mathcal O_S$ as a sub-bundle for each $k\in I_i$ , we have $(\cap_{k\in I_i}V^k)\otimes \mathcal O_S$ is the kernel of $\mathcal V_S\rightarrow \oplus_{k\in I_i}\mathcal V_S/\mathcal W_k$, thus for each $k\in I_i$, the zero morphism $(\cap_{k\in I_i}V^k)\otimes\mathcal O_S\hookrightarrow \mathcal V_S\rightarrow \mathcal V_S/\mathcal W_k $ will factor through $\mathcal W_k$, which is still a sub-vector bundle of $\mathcal V_S$. Therefore, $\cap_{k\in I_i}V^k\otimes \mathcal O_S$ is a sub-bundle of $\mathcal W_k$.

 Recall that $v_{[n]}\in (\bigcap_{i=1}^n V^i)\setminus (\bigcup_{i=1}^n V_i)$, thus $v_{[n]}\in \bigcap_{k\in I_i}V^k$. Hence $\spa(v_{[n]})\otimes\mathcal O_S$ is a sub-bundle of $\mathcal W_i$.  Then we have $\spa (V_i, v_{[n]})\otimes \mathcal O_S$ is a sub-bundle of $\mathcal W_i$, but notice that
$$\rk (\spa(V_i, v_{[n]})\otimes \mathcal O_S)=d_i+1=\rk \mathcal W_i,
$$
by Lemma \ref{vector bundle}, $\spa (V_i, v_{[n]})\otimes \mathcal O_S=\mathcal W_i$.

 By Case (i) and (ii), $\mathcal W_i=\spa(V_i, v_{[n]})\otimes\mathcal O_S$ for each $1\leq i\leq n$. We conclude that for each scheme $S/k$,
 the set $\mathfrak H(S)=\{(\spa(V_1, v_{[n]})\otimes \mathcal O_S, \cdots,\spa(V_n, v_{[n]})\otimes \mathcal O_S)\}$, a one-point functor.
\hfill  $\Box$

\begin{Cor}\label{F_1 cohen} The moduli functor $\mathfrak F_1$ is multiplicity-free and hence Cohen-Macaulay.
\end{Cor}

\noindent\emph{\textbf{Proof:}}

By Proposition \ref{vector translation} and Lemma \ref{if and only if} ,whenever $\sum_{i\in [n]}c_i=p$, we an find $\{V^i\subset V\}_{i=1}^n$, such that $\Pi_{i=1}^n\mathbb P(V^i/V_i)$ intersects with $\mathfrak F_1$ at most one point, and by Proposition \ref{answer multiplicity}, if it intersects at one point, the multiplicity will be one. Therefore we conclude the closure of the image of the rational map $\mathbb P(V)\dashrightarrow \Pi_{i=1}^n\mathbb P(V/V_i)$ is multiplicity-free, then by a result of Brion in \cite{mulfree}, it is Cohen-Macaulay.
\hfill $\Box$

\vspace{.3cm}

\section{Initial Degeneration of $X(V_1,\cdots, V_n)$}\label{initial ideal}

In this section, we will study the initial degeneration of $X(V_1,\cdots, V_n)$ under a prescribed term order.

\subsection{Genericity Setup}

Let $\{V_i\}_{i=1}^n$ be a finite sequence of subspaces of $V$, and denote $d_I=\dim (\cap_{i\in I}V_i)$.

Consider the natural map $V\rightarrow \Pi_{i=1}^n V/V_i$.  Choose coordinate system for $V$ and $\{V/V_i\}_{i=1}^n$ in such a generic way that for $A_i:\ V\rightarrow V/V_i$,  following matrix
$$A_I=\begin{bmatrix}
A_1\\
\vdots\\
A_n
\end{bmatrix}
$$
has the following property:

Let $C$ be any set of row vectors in $A_I$, which takes $a_m$ rows taken in $A_m$. Let $I_C$ to be the set $\{i|i\in [n], C\cap \{\te{row vectors in}\  A_i\}\neq \emptyset\}$.
 If for any $ I'\subset I_C$, $\sum_{m\in I'}a_m\leq r+1-d_{I'}$  ( we will call this condition as $\mathrm{Condition (\star)}$), then $\dim (\spa C)=|C|$.

\begin{rmk} Such a way of choosing basis to satisfy the genericity set-up always exists. A priori, for arbitrary choice of basis of $V$ and $V/V_i$'s whose corresponding matrix is $A'_I$,  there exists one set $C$ of row vectors with condition $(\star)$ has $\dim(\spa C)=|C|$, and row operations which keep this property for $C$ form a dense open subset of $\Pi_{i=1}^n Gl (V/V_i)$ and  denote this open subset as $U_C$.  For another set $C'$ with condition $(\star)$ and $|C'|=|C|$, we can do elementary row operation to permute the rows between $C'$ and $C$ so that $\dim(\spa C')=|C'|$, thus there is a nonempty dense open row operations in $\Pi_{i=1}^n Gl (V/V_i)$ which makes $\dim(\spa C')=|C'|$ and we will denote this subset as $U_{C'}$.  There are only finitely many $C'$ with condition $(\star)$ in the collection of subsets of the set of row vectors and for each $C'$.

Let $I_{\star}=\{\mathrm{row\ vectors\ of\ A_I\ with\ Property}\ (\star)\}$, and it is clear that $I_{\star}$ is finite. $\bigcap_{C'\in I_{\star}}U_{C'}$ is nonempty because each $U_{C'}$ is nonempty and dense open. Take $\beta\in \bigcap_{C'\in I_{\star}}U_{C'}$, then $\beta\cdot A'_I$ will satisfy the genericity setup.
\end{rmk}

\subsection{Prime Decomposition of the Initial Ideal}

Recall that $M(h)=\{(m_1,\cdots,m_n)\in \mathbb Z^n_{\geq0}:\ \forall I\subset [n],\ r+1-\sum_{i\in I}m_i> d_I\}$ and
$M(\dim(\overline{\mathrm{Im}f}))=M(p)$, where $p=\max\{h: M(h)\neq\emptyset\}$.
\begin{Def} $\widehat{M}=\{(m_1,\cdots, m_n)\in \bigcup_{h=1}^p M(h), \forall k\in [n], \exists I_k\subset [n], k\in I_k, r+1-\sum_{i\in I_k}m_i=d_{I_k}+1\}$.
\end{Def}

A priori,  $M(p)\subset \widehat M$. Since $d_I$ can be realised as the dimension of intersection of linear subspaces of a vector space, $M(p)$ and $\widehat M$ are actually \emph{equal}, which is proved in the following proposition.

\begin{Prop}\label{matroid} $M(p)=\widehat M$
\end{Prop}

\noindent\emph{\textbf{Proof:}}  Given an $(m_1,\cdots, m_n)\in\widehat M$, for each $1\leq h\leq n$,let $S_{m_h}:=\{I\subset [n]: h\in I, r+1-\sum_{i\in I}m_i=d_I+1 \}$, which is non empty by the definition of $\widehat{M}$ and forms a partially ordered set under inclusion. I claim that $S_{m_h}$ has a unique maximal element.

Suppose that $I_1$ and $I_2$ are both maximal. Consider the intersection $I_1\cap I_2$ which still contains $h$.
$$\begin{array}{rcl}
 r+1-\sum_{i\in I_1\bigcup I_2}m_i&=&(r+1-\sum_{i\in I_1}m_i)+(r+1-\sum_{i\in I_2}m_i)-(r+1-\sum_{i\in (I_1\bigcap I_2)}m_i)\\
                                  &=&(d_{I_1}+1)+(d_{I_2}+1)-(r+1-\sum_{i\in (I_1\bigcap I_2)}m_i).

 \end{array}
$$

On the other hand, $d_{I_1}=\dim(\bigcap_{i\in I_1}V_i)$ and $d_{I_2}=\dim(\bigcap_{i\in I_2}V_i)$, and both $\bigcap_{i\in I_1}V_i$ and $\bigcap_{i\in I_2}V_i$ are linear subspaces of $\bigcap_{i\in (I_1\bigcap I_2)}V_i$, thus
$$ \spa (\cap_{i\in I_1}V_i,\ \cap_{i\in I_2}V_i)\subset \bigcap_{i\in ( I_1\bigcap I_2)}V_i.
$$

$$\begin{array}{rcl}
          d_{I_1}+d_{I_2}-d_{I_1\bigcup I_2}&=&\dim(\cap_{i\in I_1}V_i)+\dim(\cap_{i\in I_2}V_i))-\dim(\cap_{i\in I_1\bigcup I_2}V_i)\\
                                            &=&\dim \spa (\cap_{i\in I_1}V_i,\ \cap_{i\in I_2}V_i)\\
                                            &\leq&\dim (\cap_{i\in ( I_1\bigcap I_2)}V_i)=d_{I_1\bigcap I_2}\\
                                            &<&(r+1-\sum_{i\in (I_1\bigcap I_2)}m_i).
\end{array}
$$
Therefore we have
$$\begin{array}{rcl}
r+1-\sum_{i\in I_1\bigcup I_2}m_i&=&(d_{I_1}+1)+(d_{I_2}+1)-(r+1-\sum_{i\in (I_1\bigcap I_2)}m_i)\\
                                 &=&(d_{I_1}+d_{I_2}-d_{I_1\bigcup I_2})+2+d_{I_1\bigcup I_2}-(r+1-\sum_{i\in (I_1\bigcap I_2)}m_i)\\
                                 &<&(r+1-\sum_{i\in (I_1\bigcap I_2)}m_i)+2+d_{I_1\bigcup I_2}-(r+1-\sum_{i\in (I_1\bigcap I_2)}m_i)\\
                                 &=&2+d_{I_1\bigcup I_2}.
\end{array}
$$
Note that $r+1-\sum_{i\in I_1\bigcup I_2}m_i>d_{I_1\bigcup I_2}$, and $r+1-\sum_{i\in I_1\bigcup I_2}m_i$ is an integer, $r+1-\sum_{i\in I_1\bigcup I_2}m_i$ has to be equal to $d_{I_1\bigcup I_2}+1$, therefore, $I_1\bigcup I_2\in S_{m_h}$, thus $I_1\subset I_1\cup I_2$, and $I_2\subset I_1\cup I_2$, but $I_1$ and $I_2$ are maximal, then $I_1=I_1\cup I_2=I_2$.   Therefore $S_{m_h}$ has a unique maximal element.

Given $(m_1,\cdots, m_n)\in\widehat{M}$, and let $I_{m_h}$ be the unique maximal element of $S_{m_h}$, then I claim (i): $I_{m_h}=I_{m_j}$ for any $j\in I_{m_h}$; (ii): $I_{m_h}\cap I_{m_j}=\emptyset$ for $j\notin I_{m_h}$.

Proof of Claim (i): For any $j\in I_{m_h}$, we have $I_{m_h}\in S_{m_j}$, and since $I_{m_j}$ is the maximal element in $S_{m_j}$, we have $I_{m_h}\subset I_{m_j}$.  Therefore, $h\in I_{m_j}$, then $I_{m_j}\in S_{m_h}$ and since $I_{m_h}$ is the maximal in $S_{m_h}$, we have $I_{m_j}\subset I_{m_h}$. We conclude that $I_{m_h}=I_{m_j}$.

Proof of Claim (ii): For any $j\notin I_{m_h}$, consider $I_{m_j}$, the maximal element of $S_{m_j}$.  If $I_{m_j}\cap I_{m_h}\neq \emptyset$, take $q\in I_{m_j}\cap I_{m_h}$, then by claim one, $I_{m_h}=I_{m_q}=I_{m_j}\ni j$, thus $j\in I_{m_h}$, a contradiction. Therefore $I_{m_h}\cap I_{m_j}=\emptyset$ for $j\notin I_{m_h}$.

By Claim (i) and (ii),given any $(m_1,\cdots,m_n)\in\widehat{M}$,  one can decompose $[n]$ into disjoint union  $I_{u_1}\sqcup I_{u_2}\cdots \sqcup I_{u_\ell}$ where $\{u_1,\cdots, u_\ell\}\subset [n]$, such that $r+1-\sum_{i\in I_{u_q}}m_i=d_{I_{u_q}}+1$ for each $1\leq q\leq \ell$.

For each $(a_1,\cdots,a_n)\in M(p)$, we have $r+1-\sum_{i\in I_{u_q}}a_i>d_{I_{u_q}}$ but $r+1-\sum_{i\in I_{u_q}}m_i=d_{I_{u_q}}+1$, thus $\sum_{i\in I_{u_q}}a_i\leq\sum_{i\in I_{u_q}}m_i$ for each $1\leq q\leq \ell$.
$$\begin{array}{rcl}
 p=\sum_{i=1}^n a_i&=&\sum_{q=1}^{\ell}\sum_{i\in I_{u_q}}a_i\\
                   &\leq&\sum_{q=1}^{\ell}\sum_{i\in I_{u_q}}m_i\\
                   &=&\sum_{i=1}^n m_i
 \end{array}.
$$
Recall that by definition, $\widehat{M}\subset \bigcup_{h=1}^p M(h)$, and therefore $\sum_{i=1}^n m_i\leq p$, thus $\sum_{i=1}^n m_i=p$. Hence for each $(m_1,\cdots,m_n)\in \widehat{M}$, we have $\sum_{i=1}^n m_i=p$, then $\widehat{M}\subset M(p)$.  Thus the proposition follows.
\hfill $\Box$
\begin{rmk} A priori, $M(p)\subset \widehat M$, and in general, this containment is strict. The proof of Proposition \ref{matroid} assumes that $d_I$ can be realized as the dimension of intersections of linear subspaces of a single vector space, but actually, by exactly the same argument, one can prove the same result for general matroid rank functions.
\end{rmk}

 \begin{conve}
 Recall $d_i=\dim V_i$. Let $k[x_{i,1}, x_{i,2}, \cdots, x_{i, r+1-d_i}]$ denote the projective coordinate ring of $\mathbb P(V/V_i)$, and our lex order is taken to be $$\left\{
 \begin{array}{cl}
x_{i,j}\succ x_{\ell,m} &\text{if}\  i>\ell;\\
x_{i,j}\succ x_{\ell,m}& \text{if}\  i=\ell,\ j>m
\end{array}
\right.
$$
\end{conve}

\begin{Def}$P_{(m_1,\cdots, m_n)}$ is a prime ideal associated to $(m_1,\cdots, m_n)$ defined as the prime ideal $$<x_{1,1}, x_{1,2}, \cdots, x_{1,r-d_1-m_1 },\ \dots , \ x_{i,1},\cdots, x_{i, r-d_i-m_i}, \ \dots,\  x_{n,1}, \cdots,x_{n, r-d_n-m_n} >$$
\end{Def}

\begin{Def}\label{I_o}  Define $I_o$ to be $\bigcap_{(m_1,\cdots,m_n)\in  M(p)}P_{(m_1,\cdots,m_2)}$.
\end{Def}
\vspace{0.3cm}
Let's first make a explicit description of $I_o$  by writing down its generators.
\begin{conve}For a monomial $x_{\ell, \ell_1}x_{2, \ell_2}\cdots x_{i, \ell_i}\cdots x_{n, \ell_n}$ in the coordinate ring of $\Pi_{1}^n \mathbb P(V/V_i)$, when $\ell_i=0$, we set up $x_{i, \ell_i}=1$.
\end{conve}
The following proposition writes down a set of generators of $I_o$.
\begin{Prop}\label{generators} $I_o$ can be generated by monomials $x_{1, \ell_1}x_{2, \ell_2}\cdots x_{i, \ell_i}\cdots x_{n, \ell_n}$ where $\ell_i$'s are non negative integers satisfying the following condition:
$$\exists I\subset [n]\  \text{such that}\ r+1-\sum_{i\in I}(r+1-d_i-\ell_i)\leq d_I
$$
\end{Prop}
Before proving this proposition, let's first prove the following lemma:
\begin{Lemma}\label{whyout}Given a monomial $x_{1, \ell_1}x_{2, \ell_2}\cdots x_{i, \ell_i}\cdots x_{n, \ell_n}$, such that for any $I\in [n]$, we have $r+1-\sum_{i\in I} (r+1-d_i-\ell_i) >d_I$, then this monomial is not in $I_o$.

\end{Lemma}
\noindent \emph{\textbf{Proof:}}
For any $I\subset [n]$, we have $r+1-\sum_{i\in I} (r+1-d_i-\ell_i) >d_I$, therefore $(r+1-d_1-\ell_1,\cdots, r+1-d_n-\ell_n)\in M(n)$ for some $n\leq p$. Hence there exists $(m_1,\cdots,m_n)\in \widehat{M}$ such that for each $1\leq i\leq n$, we have $r+1-d_i-\ell_i\leq m_i$, which implies $r-d_i-m_i<\ell_i$ for each $i\in [n]$.

Consider the prime ideal $P_{(m_1,\cdots, m_n)}$, and we conclude that this monomial  $x_{1, \ell_1}x_{2, \ell_2}\cdots x_{i, \ell_i}\cdots x_{n, \ell_n}$ is not in this prime ideal $P_{(m_1,\cdots, m_n)}$ because $x_{i, \ell_i}\notin P_{(m_1,\cdots, m_n)}$ for all $i$'s.

By Lemma \ref{matroid},$\widehat{M}=M(p)$, thus $(m_1,\cdots, m_n)\in M(p)$. The prime ideal $P_{(m_1,\cdots, m_n)}$ is a prime factor of $I_o$, therefore $x_{1, \ell_1}x_{2, \ell_2}\cdots x_{i, \ell_i}\cdots x_{n, \ell_n}$ is not in $I_o$, which is contained in $P_{(m_1,\cdots, m_n)}$.

\hfill $\Box$

\begin{Lemma}\label{whyin}Given a monomial $x_{1, \ell_1}x_{2, \ell_2}\cdots x_{i, \ell_i}\cdots x_{n, \ell_n}$, with an $I\subset [n]$ such that $$r+1-\sum_{i\in I}(r+1-d_i-\ell_i)\leq d_I,$$  then this monomial $x_{1, \ell_1}x_{2, \ell_2}\cdots x_{i, \ell_i}\cdots x_{n, \ell_n}$ is in $I_o$.
\end{Lemma}
\noindent \emph{\textbf{Proof:}}\ For any $(m_1,\cdots, m_n)\in M(p)$, we have for the $I\in [n]$ in this lemma, $r+1-\sum_{i\in I}m_i> d_I$. We also have $r+1-\sum_{i\in I}(r+1-d_i-\ell_i)\leq d_I$, so there exists $i'\in I$, such that $m_{i'}< r+1-d_{i'}-\ell_{i'}$,thus $\ell_{i'}\leq r-d_{i'}-m_{i'}$.
Therefore, $x_{i',\ell_{i'}}\in P_{(m_1,\cdots, m_n)}$, which implies $x_{1, \ell_1}x_{2, \ell_2}\cdots x_{i, \ell_i}\cdots x_{n, \ell_n}\in  P_{(m_1,\cdots, m_n)}$. Since $(m_1,\cdots, m_n)$ is an arbitrary element in $M(p)$, thus we have $x_{1, \ell_1}x_{2, \ell_2}\cdots x_{i, \ell_i}\cdots x_{n, \ell_n}\in \cap _{(m_1,\cdots, m_n)\in M(p)}P_{(m_1,\cdots, m_n)}$
\hfill  $\Box$
\begin{Lemma}\label{form} Generaters of $I_o$ can be written as monomials of the form $x_{1, \ell_1}x_{2, \ell_2}\cdots x_{i, \ell_i}\cdots x_{n, \ell_n}$
\end{Lemma}
\noindent \emph{\textbf{Proof:}}\ Intersections of monomial ideals are still monomial, thus $I_o$ is a monomial ideal.  Let $I_k$ be a subset of $\{0,1,2,\cdots, r+1-d_i\}$ and $i_k$ be the smallest integer of $I_k$ . If $\Pi_{i=1}^n\Pi_{h\in I_k}x_{i,h}\in I_o$, then $\Pi_{k=1}^n x_{k, i_k}$ is also in $I_o$.
\hfill  $\Box$

Now let's prove the Proposition \ref{generators}:

\noindent \emph{\textbf{Proof:}} \ By lemma \ref{form}, we know that $I_o$ is generated by monomials of the form $x_{1, \ell_1}x_{2, \ell_2}\cdots x_{i, \ell_i}\cdots x_{n, \ell_n}$. And by Lemma \ref{whyout} and Lemma \ref{whyin}, we know that generators of $I_o$ are of the form $x_{1, \ell_1}x_{2, \ell_2}\cdots x_{i, \ell_i}\cdots x_{n, \ell_n}$ with $r+1-\sum_{i\in I}(r+1-d_i-\ell_i)\leq d_I$ for some $I\in [n]$.
\hfill  $\Box$
\begin{Def}A generator of the form $x_{1, \ell_1}x_{2, \ell_2}\cdots x_{i, \ell_i}\cdots x_{n, \ell_n}$ for $I_o$ is called a  \emph{\textbf{irredundant generator}} if $I=\{\ell_i|\ell_i\neq 0\}$ satisfies the condition that $r+1-\sum_{i\in I}(r+1-d_i-\ell_i)\leq d_I$ and for any $I'\subsetneqq I$, $r+1-\sum_{i\in I'}(r+1-d_i-\ell_i)> d_{I'}$
\end{Def}
\begin{Lemma}\label{c1}Irredundant generators exist and all of the irredundant generators generate $I_o$.
\end{Lemma}

\noindent{\emph{\textbf{Proof:}}}
For any $x_{1, \ell_1}x_{2, \ell_2}\cdots x_{i, \ell_i}\cdots x_{n, \ell_n}\in I_o$, there is $I_1\subset [n]$, with $r+1-\sum_{i\in I_1}(r+1-d_i-\ell_i)\leq d_{I_1}$, then $\Pi_{i\in I_1}x_{i,\ell_i}\in I_o$. Take $I_2\subsetneqq I_1$ with $r+1-\sum_{i\in I_2}(r+1-d_i-\ell_i)\leq d_{I_2}$, we can truncate this generator to the form $\Pi_{i\in I_2}x_{i, \ell_i}$. Repeating this procedure, we will get :
$$ I_1\supsetneqq I_2\cdots\supsetneqq I_h\cdots\ \ \ .
$$
The sequence will stop eventually, and denote the last term in the sequence as $I_k$ then $\Pi_{i\in I_k}x_{i,\ell_i}$ is an irredundant generator. The same operation will reduce very monomial generators of the form $\Pi_{i=1}^n x_{i,\ell_i}$ into an irredundant generator.
\hfill $\Box$

In computer vision \cite{computer}, the length of the generators (their generators corresponds to ``irredundant generators'' in our language) was bounded by $4$ in Theorem 2.1. The next lemma gives a generalisation of this upper bound.
\begin{Cor}  The length of the irredundant generators of $I_o$ is bounded by $$\min\{r+1, n\}.$$
\end{Cor}
\noindent \emph{\textbf{Proof:}}\  By Lemma \ref{c1}, generators of the $I_o$ are monomials $x_{1, \ell_1}x_{2, \ell_2}\cdots x_{i, \ell_i}\cdots x_{n, \ell_n}$ with $\exists I\in [n]$ such that $r+1-\sum_{i\in I}(r+1-d_i-l_i)\leq d_I$ and for any $I'\subsetneqq I$, $r+1-\sum_{i\in I'}(r+1-d_i-\ell_i)> d_{I'}$.\\
Pick $I'\subset I$ with $|I'|=|I|-1$. By $r+1-\sum_{i\in I'}(r+1-d_i-\ell_i)> d_{I'}$ and $(r+1-d_i-\ell_i)\geq 1$ for each $i$, $|I'|\leq\sum_{i\in I'}(r+1-d_i-\ell_i)< r+1-d_{I'}\leq r+1 $, therefore $|I|=|I'|+1\leq r+1$. On the other hand, $|I|\leq n$, thus this corollary follows.
\hfill  $\Box$

Let $A_i$ be the matrix representation of the linear map $V\rightarrow V/V_i$, and $q_i$ be $(x_{i,1},\cdots, x_{i, r+1-d_i})^T$, and let $I=\{\delta_1, \cdots, \delta_{|I|}\}$ be any subset $[n]$,  and consider the following matrix associated to $I\subset [n]$,
$$\begin{bmatrix}
  A_{\delta_1}&q_{\delta_1}&0&\hdots&0 \\
  A_{\delta_2}&0&q_{\delta_2}&\vdots&0\\
  \vdots&\vdots& &&&\\
  A_{\delta_{|I|}}&0&\hdots&\hdots&q_{\delta_{|I|}}
\end{bmatrix}$$
Consider all of the $r+1-d_I+|I|$-minor of the above matrix, and denote the ideal generated by these minors as $I_f$.
\begin{Lemma}$I_o\subset\Int_{\prec} I_f$
\end{Lemma}
\noindent \emph{\textbf{Proof:}}\  It suffices to check that all of the irredundant generators of $I_o$ are in $\Int_{\prec} I_f$.
We are going to check monomials of the form $x_{1,\ell_{1}}\cdots x_{n,\ell_{y_n}}$ whose $I=\{\ell_i:\ \ell_i\neq 0\}$ has the property that $r+1-\sum_{i\in I}(r+1-d_i-\ell_i)\leq d_I$ and for any $I'\subsetneqq I,\ r+1-\sum_{i\in I'}(r+1-d_i-\ell_i)>d_{I'}$. Now let $I=\{y_1,\cdots, y_{|I|}\}$, then we need to prove that $x_{y_1, \ell_{y_1}}\cdots x_{y_{|I|},\ell_{y_{|I|}}}\in \Int_{\prec}I_f$.
Let us to consider matrix with the following form:\\
Denote  $A_{i,j}$ be the $j-$th row of $A_i$; and let $I=\{y_1,\cdots, y_{|I|}\}$\\
$$A'=\begin{bmatrix}
A_{y_1,l_{y_1}}&x_{y_1, \ell_{y_1}}&0\\
\vdots&\vdots&\vdots&0\\
A_{y_1, r+1-d_{y_1}}&x_{y_1, r+1-d_{y_1}}&0\\
A_{y_k,l_{y_k}}&0&\hdots&x_{y_k, \ell_{y_k}}&0\\
\vdots&\vdots&\vdots&\vdots\\
A_{y_k, r+1-d_{y_k}}&0&\hdots&x_{y_k, r+1-d_{y_k}}&0\\
A_{y_{|I|},l_{y_{|I|}}}&0&\hdots&0&x_{y_{|I|}, \ell_{y_{|I|}}}\\
\vdots&\vdots&\vdots&0\\
A_{y_{|I|}, r+1-d_{y_{|I|}}}&0&\hdots&0&x_{y_{|I|}, r+1-d_{y_{|I|}}}\\
\end{bmatrix}$$
In the matrix above, there are $\sum_{i\in I}(r+2-d_i-l_i)$ rows and $(r+1+|I|)$ columns, and by
$$r+1-d_I\leq \sum_{i\in I}(r+1-d_i-l_i),$$
 we have
 $$(r+1-d_I+|I|)\leq \sum_{i\in I}(r+2-d_i-l_i).$$

Then for each submatrix $A'_i$ of $A'$ where $A'_i$ is the following :
$$A'_i=\begin{bmatrix}
A_{y_i,l_{y_i}}&0&\hdots&x_{y_i, \ell_{y_i}}&0\\
\vdots&\vdots&\vdots&\vdots\\
A_{y_i, r+1-d_{y_i}}&0&\hdots&x_{y_i, r+1-d_{y_i}}&0\\
\end{bmatrix},$$
Since $r+1-\sum_{i\in I}(r+1-d_i-\ell_i)\leq d_I$, there exists $\{0\leq t_i\leq r+1-d_i-\ell_i\}_{i\in I}$ with $\sum_{i=1}^{|I|}t_i=[\sum_{i\in I}(r+1-d_i-\ell_i)]-(r+1-d_I)\geq 0$,
$$A_{I,2}=\begin{bmatrix}
A_{y_1,\ell_{y_1}}&x_{y_1, \ell_{y_1}}&0\\
\vdots&\vdots&\vdots&0\\
A_{y_1, r+1-d_{y_1}-t_{y_1}}&x_{y_1, r+1-d_{y_1}-t_{y_1}}&0\\
A_{y_k,\ell_{y_k}}&0&\hdots&x_{y_k, \ell_{y_k}}&0\\
\vdots&\vdots&\vdots&\vdots\\
A_{y_k, r+1-d_{y_k}-t_{y_k}}&0&\hdots&x_{y_k, r+1-d_{y_k}-t_{y_k}}&0\\
A_{y_{|I|},\ell_{y_{|I|}}}&0&\hdots&0&x_{y_{|I|}, \ell_{y_{|I|}}}\\
\vdots&\vdots&\vdots&0\\
A_{y_i, r+1-d_{y_{|I|}}-t_{y_{|I|}}}&0&\hdots&0&x_{y_{|I|},r+1-d_{y_{|I|}}-t_{y_{|I|}}}\\
\end{bmatrix},$$
whose number of rows is $[\sum_{i\in I}(r+2-d_i-\ell_i)]-[(\sum_{i\in I}(r+1-d_i-\ell_i))-(r+1-d_I)]=r+1-d_I+|I|$.

Consider the following matrix:

$$A_{I,3}=\begin{bmatrix}
A_{y_1,\ell_{y_1}+1}\\
\vdots\\
A_{y_1, r+1-d_{y_1}-t_{y_1}}\\
A_{y_k,\ell_{y_k}+1}\\
\vdots\\
A_{y_k, r+1-d_{y_k}-t_{y_k}}\\
A_{y_{|I|},\ell_{y_{|I|}}+1}\\
\vdots\\
A_{y_{|I|}, r+1-d_{y_{|I|}}-t_{y_{|I|}}}\\
\end{bmatrix},$$

One can see that $A_{I,3}$ has $r+1-d_I$ rows and also by the irreduncy condition and genericity setup, $\rk A_{I,3}=r+1-d_I$. Also note that $A_{I,3}$ has $(r+1)$ columns so there exists a $(r+1-d_I)\times (r+1-d_I)$ submatrix of $A_{I,3}$ with full rank and denote this matrix as $A_{I,4}$.

 Now we take $(r+1-d_I+|I|)\times (r+1-d_I+|I|)$ submatrix of $A'$ such that the coefficient of $\Pi_{y_i\in I}x_{i, y_i}$ (which is a irredundent generator for $I_p$) is $\det A_{I,4}$ which is not zero. Thus the lex initial term is $\Pi_{y_i\in I}x_{i, y_i}$.  Thus the lemma follows.

\hfill  $\Box$

\begin{Lemma}\label{prime decomposition} $\Int_{\prec} I_f=I_o$
\end{Lemma}
\emph{\textbf{Proof:}}
The closed subscheme $X_{\Int I_f}$ cut by $\Int_{\prec} I_f$ is a flat degeneration of the subscheme $X(I_f)$, thus the multi-degrees of $X(\Int_{\prec} I_f)$ and $X(I_f)$ are the same.  Also note that $X(I_o)$ and $X(I_f)$ have the same multi-degrees, therefore, $X(\Int_{\prec} I_f)$ and $X(I_o)$ have the same multi-degrees.

Since $I_o\subset \Int_{\prec} I_f$, thus $X(\Int_{\prec} I_f)\subset X(I_o)$, and $X(I_o)$ is reduced and equidimensional. Therefore $X(\Int_{\prec} I_f)=X_{I_o}$. Thus the lemma follows.
\hfill     $\Box$

\begin{Cor} The closure of the image of the rational map represents the functors $\mathfrak F_1$ and $\mathfrak F_2$. In particular, it is Cohen-Macaulay.
\end{Cor}

\noindent \emph{\textbf{Proof:}}\  Recall that the closed subscheme in Description II is cut by $I_f$, which represents $\mathfrak F_2$.
Since by Lemma \ref{prime decomposition} we have $\Int_{\prec}I_f=I_o$, which is square free, the variety cut by $I_f$ is reduced.  The closed points of the scheme of Description II is the same as the closed point of the moduli functor $\mathfrak F_1$, whose closed points are proven to be the closure of the image of the rational map in Theorem \ref{closed points F_1}. The closure of the image of the rational map is also reduced, and therefore is identical with the variety cut by $I_f$, which represent the isomorphic moduli functors $\mathfrak F_1$ and $\mathfrak F_2$.

In Corollary  \ref{F_1 cohen}, it shows that $\mathfrak F_1$ is Cohen-Macaulay, then so is the closure of the image of the rational map.

\hfill   $\Box$

\section{ hilbert polynomial of $X(V_1,\cdots, V_n)$}\label{hilbert polynomial}

It is a natural question to ask  after knowing the multi-degrees, what the Hilbert-Polynomial is. In this section, we are going to compute the Hilbert polynomial instead of  counting monomials but using degeneration techniques to give a very intuitive argument.

 Recall the Hilbert polynomial of product of projective spaces $\Pi_{i=1}^n\mathbb P^{\ell_i}$ is equal to $\Pi_{i=1}^n\binom{u_i+\ell_i}{\ell_i}$ where $u_i$'s are variables.

 And for $A$ and $B$ two closed reduced subscheme of the same product of projective spaces, let $HP(\mathcal X)$ be the Hilbert polynomial of $\mathcal X$, a reduced subscheme of a product of projective spaces, then we have $HP(A\cup B)=HP(A)+HP(B)-HP(A\cap B)$. What about for $HP(A\cup B\cup C)$?
 Plugging in the same formula, we have $HP(A\cup B\cup C)=HP(A\cup B)+HP(C)-HP((A\cup B)\cap C)$. In general we shouldn't expect that $HP((A\cup B)\cap C)=HP((A\cap C)\cup (B\cap C))$. In our case, each variety is cut out by a monomial ideal whose generators a degree one monomials, for which we have the following easy lemma:

 \vspace{.3cm}

 \begin{Lemma}\label{monomial property} Let $I_1,\cdots, I_n$ be degree one monomial ideals of a multi-variable polynomial ring. Then
 $$(I_1\cap I_2\cap\cdots\cap I_{n-1})+I_n=(I_1+I_n)\cap I_2\cdots\cap I_{n-1}+I_n
 $$
 \end{Lemma}

\noindent\emph{\textbf{Proof:}}\
Let $J$ be a monomial ideal of the same polynomial ring, then I claim $(I_1+I_n)\cap J=(I_1\cap J)+(I_n\cap J)$. Let $M_J$ be a minimal set of monomial generators of $J$. And let $S_{I_1}$ and $S_{I_n}$ be the degree one monomial generators of $I_1$ and $I_n$ respectively.

Generally for a degree one monomial ideal $I^1$ and its degree one generating set $S_{I^1}$, let $S_{I^1}\star M_J$ be the following set of monomials:
$$\{ab:\ a\in S_{I^1}, b\in M_J\  \mathrm{with} \ a\nmid b\}\bigcup \{b|\ \exists a\in S_{I^1}\ \mathrm{with}\ a|b\}
$$

Therefore, it is not hard to see that $S_{I^1}\star M_J$ is a generating set of $I^1\cap J$, and also if $S_{I^1}=A\cup B$, then $S_{I^1}\star M_J=A\star M_J\cup B\star M_J$.

Since $(I_1+I_n)$ is also a degree one monomial ideal,  generating set of $(I_1+I_n)\cap J$ can be written as $S_{I_1+I_n}\star M_J$. By the same token, a generating set for $(I_1\cap J)+(I_n\cap J)$ is\\
 $(S_{I_1}\star M_J)\bigcup (S_{I_n}\star M_J)$. Also notice that $S_{I_1+I_n}=S_{I_1}\cup S_{I_2}$, thus we have
$$\begin{array}{rcl}
                    S_{I_1+I_n}\star M_J&=&(S_{I_1}\cup S_{I_2})\star M_J \\
                                        &=&(S_{I_1}\star M_J)\bigcup (S_{I_n}\star M_J)

\end{array}
$$

So a generating set $S_{I_1+I_n}\star M_J$ of $(I_1+I_n)\cap J$ is the same as a generating set \\
$(S_{I_1}\star M_J)\bigcup (S_{I_n}\star M_J)$ of $(I_1\cap J)+(I_n\cap J)$, then the claim follows.

Now let's go back to the statement of the lemma, note that $I_2\cap\cdots\cap I_{n-1}$ is also a monomial ideal, then by the claim we just proved, we have
$$(I_1+I_n)\cap I_2\cdots\cap I_{n-1}=(I_1\cap I_2\cap\cdots\cap I_{n-1})+I_n\cap I_2\cdots\cap I_{n-1},
$$

but also $I_n\cap I_2\cdots\cap I_{n-1}\subset I_n$, finally we have
$$\begin{array}{rcl}
 (I_1+I_n)\cap I_2\cdots\cap I_{n-1}+I_n&=&(I_1\cap I_2\cap\cdots\cap I_{n-1})+I_n\cap I_2\cdots\cap I_{n-1}+I_n\\
                                        &=&(I_1\cap I_2\cap\cdots\cap I_{n-1})+I_n

\end{array}
$$

Then the lemma follows.

\hfill $\Box$

\vspace{.3cm}

\begin{Cor} \label{interchange} Let $I_1,\cdots, I_n$ be degree one monomial ideals of a multi-variable polynomial ring. Then
 $$(I_1\cap I_2\cap\cdots\cap I_{n-1})+I_n=(I_1+I_n)\cap (I_2+I_n)\cdots\cap (I_{n-1}+I_n)
 $$
\end{Cor}

\noindent\emph{\textbf{Proof:}}\ By Lemma \ref{monomial property}, $(I_1\cap I_2\cap\cdots\cap I_{n-1})+I_n=(I_1+I_n)\cap I_2\cdots\cap I_{n-1}+I_n $, note that $(I_1+I_n)$ is also a degree one monomial ideal, and $(I_1+I_n)\cap I_2\cdots\cap I_{n-1}+I_n =I_2\cap (I_1+I_n)\cdots\cap I_{n-1}+I_n $, applying Lemma \ref{monomial property}, we have
$$\begin{array}{rcl}
(I_1\cap I_2\cap\cdots\cap I_{n-1})+I_n&=&(I_1+I_n)\cap I_2\cdots\cap I_{n-1}+I_n\\
                                       &=&I_2\cap (I_1+I_n)\cdots\cap I_{n-1}+I_n\\
                                       &=&(I_2+I_n)\cap (I_1+I_n)\cdots\cap I_{n-1}+I_n\\
                                       &=&(I_1+I_n)\cap (I_2+I_n)\cdots\cap I_{n-1}+I_n
\end{array}
$$

Then by the same token and repeated use of Lemma \ref{monomial property}, we finally have
$$(I_1\cap I_2\cap\cdots\cap I_{n-1})+I_n=((I_1+I_n)\cap (I_2+I_n)\cap\cdots\cap (I_{n-1}+I_n))+I_n.
$$

Also let's note that $I_n\subset ((I_1+I_n)\cap (I_2+I_n)\cap\cdots\cap (I_{n-1}+I_n))$. Thus we obtain:
$$\begin{array}{rcl}
(I_1\cap I_2\cap\cdots\cap I_{n-1})+I_n&=&((I_1+I_n)\cap (I_2+I_n)\cap\cdots\cap (I_{n-1}+I_n))+I_n\\
                                       &=&(I_1+I_n)\cap (I_2+I_n)\cdots\cap (I_{n-1}+I_n)
\end{array}
$$

\hfill    $\Box$

Finally we can get the following convenient property for us to compute the Hilbert polynomial.

\begin{Prop}\label{comhil} Let $P_i,\ 1\leq i\leq n$ be  closed subschemes of a product of projective spaces cut by degree one monomial ideals, then
$$HP((P_1\cup P_2\cup\cdots\cup P_{n-1})\bigcap P_n)=HP((P_1\cap P_n)\cup (P_2\cap P_n)\cup\cdots \cap(P_{n-1}\cap P_n))
$$
\end{Prop}

\noindent\emph{\textbf{Proof:}}\ It is a direct application of Corollary \ref{interchange}.
\hfill $\Box$

 \begin{Thm}The multi-variable Hilbert polynomial of the variety cut by $I_f$ equals
$$  \sum_{S\subset M(p)}(-1)^{|S|-1}\Pi_{i=1}^n \binom{u_i+\ell_{S,i}}{\ell_{S,i}}
$$
where $\ell_{S,i}$ is the smallest $i$-th component of all elements of $S$.
 \end{Thm}
 \noindent
 \emph{\textbf{Proof}}:\  Hilbert Polynomial is preserved under flat degeneration, so it suffices to compute the Hilbert polynomial of $\initial_{\prec}I_f$, whose prime decomposition has already been computed in Section \ref{initial ideal}, i.e.,
 $$\initial_{\prec}I_f=\bigcap_{(m_1,\cdots,m_n)\in  M(p)}P_{(m_1,\cdots, m_n)} $$

 For $(m_1,\cdots, m_n)\in M(r') $, we have
 $$P_{(m_1,\cdots, m_n)}=<x_{1,1}, x_{1,2}, \cdots, x_{1,r-d_1-m_1 },\ \dots , \ x_{i,1},\cdots, x_{i, r-d_i-m_i}, \ \dots,\  x_{n,1}, \cdots,x_{n, r-d_n-m_n} >.
 $$

 Thus the variety cut by $\Int_{\prec}I_f$ is a union of $\Pi_{i=1}^n\mathbb P^{m_i}$ with $(m_1,\cdots, m_n)\in M(r')$. By a direct computation, for $S\subset M(p)$
 $$  \bigcap_{(m_1,\cdots, m_n)\in S}\Pi_{i=1}^n \mathbb P^{m_i}=\Pi_{i=1}^n \mathbb P^{l_{S_i}},
 $$
 where $\ell_{S,i}$ is the smallest integer among the $i$-th components of integer vectors in $S$. Then this theorem will follow from the following claim:

 Given $T_1,\cdots, T_n$ subvarieties of $\Pi_{i=1}^n\mathbb P(V/V_i)$, where $T_i$'s are cut by degree one monomial ideal, then
 $$HP(\cup_{i=1}^n T_i)=\sum_{I\subset[n]}(-1)^{|I|-1}HP(\cap_{i\in I}T_i).
 $$
 When $n=2$, we can deduce this from $HP(A\cup B)=HP(A)+HP(B)-HP(A\cap B)$ for any closed subvariety $A$ and $B$.

 When $n=k$, assume the claim is true. For $n=k+1$, we have
 $$ HP(\cup_{i=1}^{k+1}T_i)=HP(T_1)+HP(\cup_{i=2}^{k+1}T_i)-HP(P_1\cap(\cup_{i=2}^{k+1}T_i))
 $$
 By Proposition \ref{comhil}, $$HP(T_1\cap(\cup_{i=2}^{k+1}T_i))=HP(\cup_{i=2}^{k+1}(T_1\cap T_i)).$$
 Since $T_1\cap T_i$ is also cut by degree one monomial ideal, we can use induction hypothesis to get:
 $$ HP(\cup_{i=2}^{k+1}(T_1\cap T_i))=\sum_{I'}(-1)^{(|I'|-1)}HP(\cap_{ I'}(T_1\cap T_i))=\sum_{ I'}(-1)^{(|I'|-1)}HP(\cap_{i\in I'}(T_i)\cap T_1),
 $$
 where $I'$ is a subset of $[k+1]\setminus\{1\}$.
 Therefore, we finally have
 $$\begin{array}{rcl}
 HP(\cup_{i=1}^{k+1}P_i)&=&HP(T_1)+HP(\cup_{i=2}^{k+1}T_i)-HP(P_1\cap(\cup_{i=2}^{k+1}T_i))\\
                        &=&HP(T_1)+HP(\cup_{i=2}^{k+1}T_i)-\sum_{I'\subset[k+1]\setminus \{1\}}(-1)^{(|I'|-1)}HP(\cap_{i\in I'}(T_i)\cap T_1)\\
                        &=&HP(T_1)+\sum_{I'\subset[k+1]\setminus \{1\}}(-1)^{(|I'|-1)}HP(\cap_{i\in I'}T_i)+\sum_{I'\subset[k+1]\setminus \{1\}}(-1)^{(|I'|+1)-1}HP(\cap_{i\in I'}(T_i)\cap T_1)\\
                        &=&\sum_{I\subset[k+1]}(-1)^{(|I|-1)}HP(\cap_{i\in I}T_i),
 \end{array}
 $$
thus the claim follows.
 \hfill  $\Box$

\bibliographystyle{plain}

\bibliography{imagesofrationalmap}

\vspace{0.3cm}
\begin{flushleft}
\textit{Binglin Li} \\
Department of Mathematics \\
University of California, Davis \\
One Shields Avenue, Davis California 95616\\

\vspace{1pt}
blnli@math.ucdavis.edu
\end{flushleft}

\end{document}